\documentclass[aap,MSNbibl,citesort,dvips]{arximspdf}

\doi{10.1214/12-AAP847} 
\volume{23}
\issue{2}
\pubyear{2013}
\firstpage{584}
\lastpage{616}

\makeatletter

\newcommand{\N}{\mathbb{N}}
\newcommand{\R}{\mathbb{R}}
\newcommand{\E}{\mathbb{E}}
\newcommand{\p}{\mathbb{P}}
\newcommand{\var}{\operatorname{Var}}
\newcommand{\cov}{\operatorname{Cov}}

\newtheorem{proposition}{Proposition}
\newtheorem{theorem}{Theorem}

\newtheorem{corollary}{Corollary}

\newtheorem{exactformula}{Exact Formula}
\newtheorem{compineq}{Comparison Inequality}

\newtheorem{supplem}{Supplementary Lemma}

\setattribute{abstract}    {skip} {20\p@}

\makeatother

\begin{document}
\begin{frontmatter}

\title{Bounds on the suprema of Gaussian processes, and omega results
for the sum of a random multiplicative function}
\runtitle{Bounds on the suprema of Gaussian processes}

\begin{aug}
\author[A]{\fnms{Adam J.} \snm{Harper}\corref{}\thanksref{t1}\ead[label=e1]{A.J.Harper@dpmms.cam.ac.uk}}
\runauthor{A. J. Harper}
\affiliation{University of Cambridge}
\address[A]{Department of Pure Mathematics\\
\quad and Mathematical Statistics\\
University of Cambridge\\
Wilberforce Road, Cambridge CB3 0WA\\
United Kingdom\\
\printead{e1}} 
\end{aug}

\thankstext{t1}{Supported by a studentship from the Engineering and
Physical Sciences Research Council of the United Kingdom.}

\received{\smonth{1} \syear{2011}}
\revised{\smonth{1} \syear{2012}}

%
\begin{abstract}
We prove new lower bounds for the upper tail probabilities of suprema
of Gaussian processes. Unlike many existing bounds, our results are not
asymptotic, but supply strong information when one is only a little
into the upper tail. We present an extended application to a Gaussian
version of a random process studied by Hal\'{a}sz. This leads to much
improved lower bound results for the sum of a random multiplicative
function. We further illustrate our methods by improving lower bounds
for some classical constants from extreme value theory, the Pickands
constants $H_{\alpha}$, as $\alpha\rightarrow0$.
\end{abstract}

%

\end{frontmatter}

\section{Introduction}\label{sec1}
Let $\mathcal{T}$ be a nonempty set, $(\Omega,\mathcal{F},\p)$ a
probability space and for each $t \in\mathcal{T}$ let $Z(t)$ be a
random variable defined on $(\Omega,\mathcal{F},\p)$. Suppose that
for any finite subset $\{t_{1}, t_{2},\ldots, t_{n}\} \subseteq\mathcal
{T}$, the random variable $(Z(t_{1}),\ldots,Z(t_{n}))$ has an $n$-variate
normal distribution. We will then say, a little loosely, that $\{Z(t)\}
_{t \in\mathcal{T}}$ is a \textit{Gaussian process} with parameter set
$\mathcal{T}$. We refer the reader to the book of Lifshits~\cite{lif}
for a general introduction to the theory of Gaussian processes.

In this paper we will be concerned with $\sup_{t \in\mathcal{T}}
Z(t)$, and in particular with giving lower bounds for the probability
that it is quite large. Results of this type have many applications and
the author's interest in them stems from a number-theoretic problem
that will be described later. For overviews of results in this area we
refer to two important books by Leadbetter, Lindgren and Rootz\'{e}n
\cite{llr} and by Piterbarg~\cite{pit}.

Suppose that $\mathcal{T}$ is a finite set, so that $\sup_{t \in
\mathcal{T}} Z(t)$ is certainly a genuine random variable, and
\[
\p\Bigl(\sup_{t \in\mathcal{T}} Z(t) > u\Bigr)
\]
is the probability that a multivariate normal random vector takes
values in a certain subset of $\R^{\#\mathcal{T}}$. We will also be
interested in processes with infinite index sets but will study these
by looking at suitably chosen finite subsets of points $t$. Unless the
mean vector and covariance matrix of $\{Z(t)\}_{t \in\mathcal{T}}$
have special forms, it is typically very difficult to compute the tail
probability exactly. Nevertheless, existing results offer two broad
options for lower bounding $\p(\sup_{t \in\mathcal{T}} Z(t) > u)$.
\begin{itemize}
\item One can use metric entropy/capacity methods, such as Sudakov's
minoration, to bound $\E\sup_{t \in\mathcal{T}} Z(t)$ (see Lifshits
\cite{lif}, Section 14). Together with suitable concentration
inequalities, such as that of Borell/Sudakov--Tsyrelson, this yields
explicit lower bounds on $\p(\sup_{t \in\mathcal{T}} Z(t) > u)$ for
fixed~$u$.

\item One can use techniques such as the method of comparison (which we
discuss more below), Pickands' method of double sums or Rice-type
methods (based on calculation of moments) to estimate the probability
asymptotically as $u \rightarrow\infty$ (see Piterbarg's book~\cite{pit}).
\end{itemize}

The methods listed can be powerful when attacking certain problems, but
have some unfortunate limitations. The lower bounds that one obtains
for $\E\sup_{t \in\mathcal{T}} Z(t)$ are typically off from the
truth by a multiplicative factor, and then the lower bounds for $\p
(\sup_{t \in\mathcal{T}} Z(t) > u)$ are very far from the truth for
moderately sized~$u$. The asymptotic techniques ultimately rely on,
among other things, the fact that as \mbox{$u \rightarrow\infty$}, any
correlations among the $Z(t)$ that are not perfect $\pm1$ correlations
have an increasingly negligible effect on the tail behavior. (Readers
familiar with, e.g., Berman's theorem should find this reasoning
familiar.) Unfortunately $u$ may need to be extremely large before the
techniques guarantee this effect to occur.

Piterbarg~\cite{pit} does not formulate the method of double sums or
the method of moments (for lower bounds) for fixed $u$, and the general
philosophy of those methods, that one need not analyze correlations of
$\{Z(t)\}_{t \in\mathcal{T}}$ except for extremely large
correlations, seems unsuited to obtaining such results. His version of
the method of comparison involves unspecified constants that appear to
depend on $\{Z(t)\}_{t \in\mathcal{T}}$, so one must wait for $u$ to
be sufficiently large, in an unspecified sense, before it comes into
play. (We present some normal comparison inequalities in Section \ref
{sec3}, and
when the author tried to study our Section~\ref{sec6} example using
them, he
could only show that the supremum there is larger than about $\log\log
x /\sqrt{2}$ with high probability, by studying points $t$ with
spacing $1/\sqrt{\log x}$. Our Corollary~\ref{co2} shows that supremum is
larger than about $\log\log x$ with high probability.)

In this paper we develop an alternative approach to lower bounding the
upper tail probability. The ingredients are an initial conditioning
step, followed by a comparison (in the sense of the method of
comparison) with a ``model'' Gaussian process that can be explicitly
analyzed. The resulting bounds are clean and can be nontrivial for
moderately sized $u$. Indeed, in our number-theoretic application we
will have nontrivial bounds for $u$ just larger than $\E\sup_{t \in
\mathcal{T}} Z(t)$ (and, in particular, will be able to identify the
expectation up to second order terms). As our bounds are completely
explicit, they also give information about the ``mysterious'' constants
in some asymptotic results, and our other application is a new lower
bound for the Pickands constants (defined later).

We begin with the following straightforward result.
%
%
\begin{proposition}[(Conditioning step)]\label{pr1}
Let $\{Z(t_{i})\}_{1 \leq i \leq n}$ be jointly multivariate normal
random variables. Set $r_{i,j}:=\E Z(t_{i})Z(t_{j})$, and suppose that:
\begin{itemize}
\item(\textit{centralization and normalization}) $\E Z(t_{i})=0$ and
$\E
Z(t_{i})^{2}=1$ for all $1 \leq i \leq n$;

\item(\textit{no repeated variables}) $|r_{i,j}|<1$ whenever $i \neq j$.
\end{itemize}
Then for any $u \geq0$ and any $H \geq0$,
\[
\p\Bigl(\max_{1 \leq i \leq n} Z(t_{i}) > u\Bigr) \geq\frac{H
e^{-(u+H)^{2}/2}}{\sqrt{2\pi}} \sum_{m=1}^{n} \inf_{0 \leq h \leq
H} P(m,h),
\]
where $P(m,h)$ is
\[
\p\biggl(V_{j} \leq\frac{u - r_{j,m}(u+h)}{\sqrt{1 - r_{j,m}^{2}}}
\ \forall j \leq m-1 \biggr),
\]
and the $V_{j}=V_{j,m}$ are centralized, normalized, jointly
multivariate normal random variables with correlations
\[
\frac{r_{j,k} - r_{j,m}r_{k,m}}{\sqrt{(1 - r_{j,m}^{2}) (1 -
r_{k,m}^{2})}}.
\]
\end{proposition}

We give the short proof of Proposition~\ref{pr1} in Section~\ref{sec2}.
The author had a
more involved proof of (a result like) Proposition~\ref{pr1}, based on a
``reversal of roles'' in the normal comparison procedure. Since we will
need some normal comparison results later, we present these in Section
\ref{sec3} and give a very brief description of the reversal of roles approach
as well.

We now turn to the problem of what we will be able to say about
$P(m,h)$. If the correlation structure of $\{Z(t_{i})\}_{1 \leq i \leq
m}$ is arbitrary, the answer may be essentially nothing, in which case
our attempt to give lower bounds will be at an end. However, under some
conditions on the correlation structure we can be more optimistic, and
to show this we formulate the following result.

%
\begin{proposition}[(Comparison step)]\label{pr2}
Let $u \geq0$, and suppose $h$ is sufficiently small that all the
upper bounds $(u - r_{j,m}(u+h))/\sqrt{1 - r_{j,m}^{2}}$ in the
definition of $P(m,h)$ are nonnegative. Suppose there exist numbers
$c_{j}=c_{j}(m,h), d_{j}=d_{j}(m,h)>0$ such that:
\begin{longlist}
\item$c_{j}/d_{j}$ is a nondecreasing sequence, $1 \leq j \leq m-1$;
\item$c_{\min\{j,k\}} d_{\max\{j,k\}}$ is a strict lower bound for
$r_{j,k}-r_{j,m}r_{k,m}$, for each pair $1 \leq j,k \leq m-1$.
\end{longlist}
Then for any $\delta\geq0$,
\[
P(m,h) \geq\int_{-B(\delta)}^{B(\delta)} \frac{e^{-t^{2}/2}}{\sqrt
{2\pi}}\,dt \cdot\prod_{j=1}^{m-1} \Phi\biggl(\frac{(1-\delta
)(u-r_{j,m}(u+h))}{\sqrt{1 - r_{j,m}^{2} - c_{j}d_{j}}} \biggr),
\]
where $B(\delta)=\delta\sqrt{\frac{d_{m-1}}{c_{m-1}}} \min_{1 \leq
j \leq m-1} \frac{u-r_{j,m}(u+h)}{d_{j}}$, and $\Phi$ denotes the
standard normal cumulative distribution function.
\end{proposition}

We will prove Proposition~\ref{pr2} in Section~\ref{sec4} by explicitly
constructing a
collection of Gaussian random variables with the lower bound
correlation structure suggested by the $c_{j},d_{j}$, and applying a
Brownian motion maximal inequality to analyze those. The reader might
think of this procedure as pulling out some of the dependence among the
$V_{j}$, to be analyzed nontrivially using the maximal inequality. By
doing this we gain the subtracted terms $c_{j}d_{j}$ in the product,
which will be very important, at the fairly small cost of introducing
the factor involving $B(\delta)$ [and the multiplier $(1-\delta)$].

The reader may wonder where the $c_{j},d_{j}$ will come from and
whether the lower bound obtained will not be hopelessly small in
situations of interest. In fact we can quickly deduce the following
from Propositions~\ref{pr1} and~\ref{pr2}.
%
%
\begin{theorem}\label{th1}
Let $\{Z(t_{i})\}_{1 \leq i \leq n}$ be as in Proposition~\ref{pr1}, and
suppose further that the sequence is \textup{stationary}, that is, that
$r_{j,k}=r(|j-k|)$ for some function $r$. Let $u \geq1$, and suppose
that:
\begin{itemize}
\item$r(m)$ is a decreasing nonnegative function;

\item$r(1)(1+2u^{-2})$ is at most $1$.
\end{itemize}
Then 
\begin{eqnarray*}
\p\Bigl(\max_{1 \leq i \leq n} Z(t_{i}) > u\Bigr)
&\geq& n \frac{e^{-u^{2}/2}}{40 u} \min\Biggl\{1,\sqrt{\frac
{1-r(1)}{u^{2}r(1)}}\Biggr\}\\
&&{}\times \prod_{j=1}^{n-1} \Phi\biggl(u\sqrt{1-r(j)}
\biggl(1+O\biggl(\frac{1}{u^{2}(1-r(j))}\biggr) \biggr) \biggr),
\end{eqnarray*}
where the implicit constant in the ``big Oh'' notation is absolute
$[$in particular, not depending on $\{Z(t_{i})\}_{1 \leq i \leq n}]$,
and could be found explicitly.\vadjust{\goodbreak}
\end{theorem}

Theorem~\ref{th1} follows by taking $H=u^{-1}$, $\delta=\min\{
u^{-2},\sqrt {r(1)/u^{2}(1-r(1))}\}$, $c_{j}=r_{j,m}=r(|m-j|)$,
$d_{j}=1-r_{j,m}=1-r(|m-j|)$ in the preceding propositions. In this
case if we did not have $c_{j}d_{j}$ in the denominators in
Proposition~\ref{pr2}, then $\sqrt{1-r(j)}$ would need to be replaced
by $\sqrt {(1-r(j))/(1+r(j))}$ in the product. We do not actually use
the theorem in this paper, as our examples require slightly different
parameter choices. However, a reader familiar with classical limit
theory for suprema of stationary processes (see, e.g., Leadbetter,
Lindgren and Rootz\'{e}n~\cite{llr}, Chapter 4) may find it instructive
to compare with those results. We may not expect to obtain precisely
sharp bounds from Theorem~\ref{th1}, because of the factor
$\min\{1,\sqrt {(1-r(1))/u^{2}r(1)}\}$, but it will supply good bounds
provided $u$ is large enough that the product term is at least $1/2$,
say. For given $r(j)$ this may be a much weaker requirement on $u$ than
in proofs of the classical results, which rely on normal comparison
inequalities. [The bound in Theorem~\ref{th1} is seen to be good
because, since we assumed that $r(m)$ is nonnegative, the tail
probability cannot be larger than $1-\Phi(u)^{n} = O(ne^{-u^{2}/2}/u)$.
An unfamiliar reader may deduce this from Comparison Inequality
\ref{comp2} in Section~\ref{sec3.1}.]

We now move on to our two examples which we hope will illustrate the
usefulness of Propositions~\ref{pr1} and~\ref{pr2}. In the theory of Gaussian
processes, much attention has been paid to (mean zero, variance one)
stationary processes whose covariance function satisfies
\[
r(t) = 1-C|t|^{\alpha} + o(|t|^{\alpha}) \qquad\mbox{as } t \rightarrow0,
\]
where $C > 0$ and $0 < \alpha\leq2$. In particular, a 1969 theorem of
Pickands~\cite{pickands1} describes the asymptotic behavior of suprema
of such processes; if $h > 0$ is fixed and if $\sup_{\varepsilon\leq t
\leq h} r(t) < 1$ for all $\varepsilon> 0$, then
\[
\lim_{u \rightarrow\infty} e^{u^{2}/2} u^{1-2/\alpha} \p\Bigl(\sup_{0
\leq t \leq h} Z(t) > u\Bigr) = \frac{hC^{1/\alpha}H_{\alpha}}{\sqrt
{2\pi}},
\]
where $H_{\alpha}$ is the so-called \textit{Pickands constant}. In a
second paper~\cite{pickands2}, Pickands used a result like this to
determine the limiting distribution, as $T \rightarrow\infty$, of a
scaled version of $\sup_{0 \leq t \leq T} Z(t)$. The
scaling\vspace*{1pt} in that theorem thus involves $H_{\alpha}$ (see,
e.g., the paper of Shao~\cite{shao} for further discussion of the role
of $H_{\alpha}$).

It appears that not very much is known about the size of $H_{\alpha}$.
Burnecki and Michna~\cite{burnmich} describe as ``mathematical
folklore'' the conjecture that $H_{\alpha}=1/\Gamma(1/\alpha)$, but
this is only known to hold for $\alpha= 1,2$. Bounds are available
more generally; for example, Shao~\cite{shao} used a representation of
$H_{\alpha}$ in terms of a nonstationary process, and various
techniques from Gaussian process theory, to show that
\begin{eqnarray*}
&&\biggl(\frac{\alpha}{4}\biggr)^{1/\alpha}
\biggl(1-e^{-1/\alpha}\biggl(1+\frac{1}{\alpha
}\biggr)\biggr)\\
&&\qquad
\leq H_{\alpha} \leq\alpha^{1/\alpha}\bigl(2.41 \sqrt{8.8 - \alpha
\log(0.4 + 2.5/\alpha)} + 0.77\sqrt{\alpha}\bigr)^{2/\alpha}
\end{eqnarray*}
when $0 < \alpha< 1$, and other bounds when $1 \leq\alpha\leq2$. D\c
{e}bicki and Kisowski~\cite{dk} subsequently improved the upper bound
on the range $1 < \alpha< 2$. D\c{e}bicki, Michna and Rolski \cite
{dmr} proved that
\[
\frac{\alpha}{8 \Gamma(1/\alpha)} \biggl(\frac{1}{4}\biggr)^{1/\alpha} \leq
H_{\alpha},\qquad 0 < \alpha\leq2,
\]
and in a 2009 preprint Michna~\cite{michna} improved this by a
multiplicative factor of~2. Note that, since $\Gamma(1/\alpha) \sim
\sqrt{2\pi\alpha} (1/e\alpha)^{1/\alpha}$ as $\alpha\rightarrow
0$, this is a much stronger bound than that of Shao~\cite{shao} under
that limit process.

Applying our methods, in Section~\ref{sec5} we improve the lower bound results
as $\alpha\rightarrow0$.
%
%
\begin{corollary}\label{co1}
There is an absolute constant $c > 0$, which could be found explicitly,
such that $H_{\alpha} \geq c\sqrt{\alpha} (e \alpha/2)^{1/\alpha}$
for all $0 < \alpha\leq2$.
\end{corollary}

For our main example, we give a detailed study of the following process:
\[
\sum_{p \leq x} g_{p} \frac{\cos(t \log p)}{p^{1/2+1/\log x}},\qquad
t \in \R,
\]
where the summation is restricted to prime numbers $p$, $g_{p}$ are
independent standard normal random variables and $x$ is a further large
parameter.

The motivation for studying this is its connection with a
number-theoretic problem of Wintner~\cite{wint}. Let $\varepsilon_{p}$
be a sequence of independent Rademacher random variables [so that $\p
(\varepsilon_{p}=1)=\p(\varepsilon_{p}=-1)=1/2$] and construct a ``random
multiplicative function'' from these, as
\[
f(n):= \cases{
\displaystyle \prod_{p|n} \varepsilon_{p}, &\quad if $n$
is squarefree,\vspace*{2pt}\cr
0, &\quad otherwise.}
\]
We also set $M(x):=\sum_{n \leq x} f(n)$. One can view $f(n)$ as a
heuristic model for some deterministic functions occurring in number
theory, such as the M\"{o}bius function. There has been quite a lot of
recent work on the behavior of $f(n)$ (e.g., by Chatterjee and
Soundararajan~\cite{chatsound}, Harper~\cite{harper}, Hough \cite
{ho2} and Lau, Tenenbaum and Wu~\cite{ltw}). However, the best known
lower bound result for $|M(x)|$ remains that of Hal\'{a}sz~\cite{hal},
who proved in 1982 that there exists a constant $B > 0$ such that,
almost surely,
\[
M(x) \neq O\bigl(\sqrt{x} e^{-B\sqrt{\log\log x \log\log\log x}}\bigr)
\qquad\mbox{as } x \rightarrow\infty.
\]
His proof, discussed in Appendix~\ref{secA}, shows that lower bound information
about the supremum of a certain Rademacher process (essentially the
process above, with the $g_{p}$ replaced by independent Rademacher
random variables) can be translated into lower bound information about $|M(x)|$.

In Section~\ref{sec6}, we use Propositions~\ref{pr1} and~\ref{pr2} to
prove results like the following.

%
\begin{corollary}\label{co2}
As $x \rightarrow\infty$,
\begin{eqnarray*}
&&\p\biggl(\sup_{1 \leq t \leq2(\log\log x)^{2}} \sum_{p \leq x} g_{p}
\frac{\cos(t \log p)}{p^{1/2+1/\log x}} \\
&&\qquad\leq\log\log x - \log\log
\log x + O((\log\log\log x)^{3/4})\biggr)
\end{eqnarray*}
is $O((\log\log\log x)^{-1/2})$.
\end{corollary}

We stress that Corollary~\ref{co2} is \textit{not} an asymptotic result
for a
single Gaussian process, but a statement about an infinite sequence of
processes depending on the parameter $x$. As $x$ grows, the variance of
the random sum grows for each fixed $t$, but also the correlation at
nearby values of $t$ decreases. (The reader may wish to look back at
these comments after he or she has read Section~\ref{sec6.1}.) For each
$x$, we
now know the supremum will typically exceed the level $\log\log x -
\log\log\log x + O((\log\log\log x)^{3/4})$; standard methods show
that the supremum is at most $\log\log x + \log\log\log x$ (say)
with probability $1-o(1)$, so Corollary~\ref{co2} is very precise in this
respect. (For each $x$ the process is ``almost'' stationary, as
explained in Section~\ref{sec6.1}, and a simple adaptation of Rice's formula
yields upper bounds for its supremum.) This precision is crucial if one
wishes to deduce things about $|M(x)|$; indeed it is the size of the
second order subtracted term $\log\log\log x$, together with the size
of the interval over which the supremum is taken, that determines what
can be said.

Together with a suitable version of the multivariate central limit
theorem, given in Appendix~\ref{secB}, Corollary~\ref{co2} allows a substantial
improvement of Hal\'{a}sz's~\cite{hal} result about $M(x)$. However,
it is possible to do better still.
%
%
\begin{corollary}\label{co3}
Let $A > 2.5$, and let $M(x)$ be the summatory function of a Rademacher
random multiplicative function, as above. It almost surely holds that
\[
M(x) \neq O\bigl(\sqrt{x}(\log\log x)^{-A}\bigr).
\]
\end{corollary}

Corollary~\ref{co2} implies Corollary~\ref{co3} with the restriction $A
> 3$, and this
is proved in Section~\ref{sec6}. The proofs in Section~\ref{sec6} are a
bit fiddly,
mostly because we must handle error terms arising in estimates for
prime number sums. A~``repeated sampling'' argument is used to deduce
Corollary~\ref{co2} in the form we need, and care is also needed to arrange
that the multivariate central limit theorem applies. However,
ultimately our results follow simply by substituting some correlation
values and parameter choices into Propositions~\ref{pr1} and~\ref{pr2}.
To prove
Corollary~\ref{co3} for all $A > 2.5$, an argument by contradiction is needed
to slightly sharpen the result of Proposition~\ref{pr2}. This is given in
Section~\ref{sec7}.\vadjust{\goodbreak}

It seems extremely likely that, almost surely, $M(x) \neq O(\sqrt
{x})$, and perhaps $M(x)$ almost surely has fluctuations of order
$\sqrt{x \log\log x}$ (by analogy with Kolmogorov's law of the
iterated logarithm). In fact, $M(x)$ might well exhibit even larger
fluctuations, since its probability distribution may have rather heavy
tails (see, e.g., Harper's article~\cite{harper}). However, an
argument like our own, ultimately based on studying a certain average
of $M(x)$ (see Appendix~\ref{secA} for justification of this comment), seems
unable to detect these large but rare fluctuations.%

We presented Proposition~\ref{pr2} in its current form, involving parameters
$c_{j},d_{j},\delta$, because this seems both easy to appreciate and
to lead to good results. However, as mentioned above, to prove the full
version of Corollary~\ref{co3} it is necessary to slightly strengthen
Proposition~\ref{pr2}. Such a strengthening may also be possible in the context
of Corollary~\ref{co1}; some of the initial steps of the Section \ref
{sec7} argument
transfer to that situation, but it is not clear whether the whole
argument goes through (except that it does not trivially do so).

The author also believes that there will be other Gaussian processes to
which Propositions~\ref{pr1} and~\ref{pr2} could usefully be applied
and hopes that the
reader might have some examples at hand.

\section{\texorpdfstring{Proof of Proposition \protect\ref{pr1}}{Proof of Proposition 1}}\label{sec2}
In view of the decomposition
\[
\p\Bigl(\max_{1 \leq i \leq n} Z(t_{i}) > u\Bigr) = \sum_{m=1}^{n} \p
\bigl(Z(t_{m})>u, Z(t_{j}) \leq u\ \forall j \leq m-1\bigr),
\]
it will suffice to show that, for any $1 \leq m \leq n$ and any $H \geq0$,
\[
\p\bigl(Z(t_{1}),\ldots,Z(t_{m-1}) \leq u, Z(t_{m}) > u\bigr) \geq\frac{H
e^{-(u+H)^{2}/2}}{\sqrt{2\pi}} \inf_{0 \leq h \leq H} P(m,h).
\]

It is well known (and easy to check, by computing correlations) that
$Z(t_{m})$ is independent of the collection of random variables
\[
Z(t_{j}) - r_{j,m}Z(t_{m}),\qquad 1 \leq j \leq m-1.
\]
These have mean zero and correlations
\[
r_{j,k}-r_{j,m}r_{k,m},\qquad 1 \leq j,k \leq m-1,
\]
and, in particular, none of them are degenerate (by assumption in
Proposition~\ref{pr1}). Thus $\p(Z(t_{1}),\ldots,Z(t_{m-1}) \leq u,
Z(t_{m}) >
u)$ is at least
\[
\int_{u}^{u+H} \p\bigl(Z(t_{j}) - r_{j,m}Z(t_{m}) \leq u - r_{j,m}x
\ \forall1 \leq j \leq m-1\bigr) \frac{e^{-x^{2}/2}}{\sqrt{2\pi}} \,dx,
\]
from which the proposition follows.\qed\vspace*{8pt}
%

In our applications, it will turn out that
\[
\p\bigl(Z(t_{1}),\ldots,Z(t_{m-1}) \leq u, Z(t_{m}) > u+H\bigr)
\]
decreases very rapidly as $H$ increases. Indeed, we will always choose
$H$ so that its effect in $P(m,H)$ is negligibly small, and therefore
only really need to understand $P(m,0)$. This is
the point of
introducing the initial decomposition of $\p(\max_{1 \leq i \leq n}
Z(t_{i}) > u)$, rather than trying to understand $\p(Z(t_{i}) \leq u
\ \forall1 \leq i \leq n)$ directly by conditioning.

\section{Normal comparison results}\label{sec3}

\subsection{Classical comparison results}\label{sec3.1}
In this subsection we present the equality underlying normal comparison
results and state some fairly classical consequences of this. We will
use these in a few places, and hopefully they will also give an
unfamiliar reader some idea of how the method of comparison, as it is
referred to by Piterbarg~\cite{pit}, is traditionally employed. Our
treatment largely follows Li and Shao~\cite{lishao}, although we would
also like to draw attention to a 1954 paper of Plackett~\cite{plac}
which contains a similar presentation of the basic comparison result.
(Plackett was interested in the numerical approximation of multivariate
normal probabilities, but some later comparison results are readily
obtained from his paper. Unfortunately this work does not seem to be
very widely known.)

If $\tilde{a}, \tilde{b} \in\R^{n}$, write $\tilde{a} \leq\tilde
{b}$ to mean that every component of $\tilde{a}$ is at most the
corresponding component of $\tilde{b}$. We have the following
identity, which is the key part of the proofs of various normal
comparison results.
%
%
\begin{exactformula}[(Following Li and Shao, and others)]\label{ex1}
Let $\tilde{X}=(X_{1},\ldots, X_{n})$ and
$\tilde{W}=(W_{1},\ldots,W_{n})$ be centralized and
normalized $n$-variate normal vectors, with covariance matrices
$\var(\tilde{X})= (\cov(X_{i},X_{j}))_{1 \leq i,j \leq n}=
(r_{i,j}^{(1)})$ and $\var(\tilde{W})= (r_{i,j}^{(0)})$ that are
nonsingular. Let $\tilde{u} \in\R^{n}$. Then
\begin{eqnarray*}
&&\p(\tilde{X} \leq\tilde{u}) - \p(\tilde{W} \leq\tilde{u})\\
&&\qquad= \sum_{1 \leq i < j \leq n} \bigl(r_{i,j}^{(1)} - r_{i,j}^{(0)}\bigr)
\int
_{0}^{1} \phi\bigl(u_{i},u_{j}; r_{i,j}^{(h)}\bigr)
\\
&&\hspace*{103pt}\qquad\quad{}\times\p\bigl(\tilde{Z}^{(h)} \leq
\tilde{u} | Z^{(h)}_{i}=u_{i}, Z^{(h)}_{j}=u_{j}\bigr) \,dh,
\end{eqnarray*}
where $\tilde{Z}^{(h)}=(Z^{(h)}_{1},\ldots,Z^{(h)}_{n})$ is multivariate
normal with covariance matrix
\[
\bigl(r_{i,j}^{(h)}\bigr):= h\var(\tilde{X}) + (1-h)\var(\tilde{W}),
\]
and $\phi(x,y; r)$ denotes the standard bivariate normal density with
correlation $r$, namely,
\[
\frac{1}{2\pi\sqrt{1-r^{2}}} e^{-(x^{2}-2rxy+y^{2})/2(1-r^{2})}.
\]
\end{exactformula}

To prove the formula one writes
\[
\p(\tilde{X} \leq\tilde{u}) - \p(\tilde{W} \leq u) = \int
_{0}^{1} \frac{d}{dh} \p\bigl(\tilde{Z}^{(h)} \leq\tilde{u}\bigr)
\,dh,
\]
observing that
\begin{eqnarray*}
\frac{d}{dh} \p\bigl(\tilde{Z}^{(h)} \leq\tilde{u}\bigr) &=&
\sum_{1 \leq i <
j \leq n} \frac{\partial}{\partial r_{i,j}^{(h)}} \p\bigl(\tilde
{Z}^{(h)} \leq\tilde{u}\bigr) \,\frac{\partial r_{i,j}^{(h)}}{\partial h} \\
&=&
\sum_{1 \leq i < j \leq n} \bigl(r_{i,j}^{(1)} - r_{i,j}^{(0)}\bigr) \int
_{-\infty}^{\tilde{u}} \frac{\partial^{2} f_{h}}{\partial y_{i}\,
\partial y_{j}} \,d\tilde{y}.
\end{eqnarray*}
Here $f_{h}$ is the density function of $\tilde{Z}^{(h)}$, the range
of integration has its obvious meaning, and the second equality uses the
fact that
\[
\frac{\partial f_{h}}{\partial r_{i,j}^{(h)}} = \frac{\partial^{2}
f_{h}}{\partial y_{i} \,\partial y_{j}},
\]
which follows by expressing the multivariate normal density in terms of
its characteristic function.

Exact Formula~\ref{ex1} provides rigorous support for the intuitive
idea that
distributions with ``nearby'' covariance matrices may have like
behavior. The inequalities that we derive next may express this in a
more striking way; they are a composite of results of Li and Shao \cite
{lishao} and of Leadbetter, Lindgren and Rootz\'{e}n~\cite{llr},
although in most respects are
%
unchanged from bounds of Slepian, Berman and Cram\'{e}r from the 1960s
(see Leadbetter, Lindgren and Rootz\'{e}n's book for the history and
references).
%
%
\begin{compineq}[(Following Leadbetter, Lindgren and Rootz\'{e}n, and Li
and Shao)]\label{comp1}
If $\tilde{X}, \tilde{W}, \tilde{u}$ are as in Exact Formula~\ref{ex1},
and $\mathbf{1}$ denotes the indicator function, then each of the
following is an upper bound for $\p(\tilde{X} \leq\tilde{u}) -
\p(\tilde{W} \leq\tilde{u})$:
\begin{longlist}
\item
\[
\frac{1}{2\pi} \sum_{1 \leq i < j \leq n} \mathbf{1}_{r_{i,j}^{(1)}
> r_{i,j}^{(0)}} \int_{r_{i,j}^{(0)}}^{r_{i,j}^{(1)}} \frac{1}{\sqrt
{1-t^{2}}} e^{-(u_{i}^{2}+u_{j}^{2})/2(1+|t|)} \,dt;
\]

\item
\begin{eqnarray*}
&&
\frac{1}{2\pi} \sum_{1 \leq i < j \leq n} \mathbf{1}_{r_{i,j}^{(1)}
> r_{i,j}^{(0)}} \bigl(\arcsin\bigl(r_{i,j}^{(1)}\bigr) -
\arcsin\bigl(r_{i,j}^{(0)}\bigr)\bigr)\\
&&\hspace*{52pt}{}\times
e^{-(u_{i}^{2}+u_{j}^{2})/2(1+\max\{|r_{i,j}^{(1)}|,|r_{i,j}^{(0)}|\}
)};
\end{eqnarray*}

\item
\begin{eqnarray*}
&&
\frac{2}{\pi} \sum_{1 \leq i < j \leq n} \mathbf{1}_{r_{i,j}^{(1)}
> r_{i,j}^{(0)}} \frac{(1+\max\{|r_{i,j}^{(1)}|,|r_{i,j}^{(0)}|\}
)^{3/2}}{(u_{i}^{2}+u_{j}^{2}) \sqrt{1-\max\{
|r_{i,j}^{(1)}|,|r_{i,j}^{(0)}|\}} }\\[-2pt]
&&\hspace*{46.2pt}{}\times e^{-(u_{i}^{2}+u_{j}^{2})/2(1+\max
\{|r_{i,j}^{(1)}|,|r_{i,j}^{(0)}|\})}.
\end{eqnarray*}
\end{longlist}
\end{compineq}

To obtain the first bound, we overestimate the conditional probability
in Exact Formula~\ref{ex1} trivially by 1 and insert the definition of
$\phi
(u_{i},u_{j}; r_{i,j}^{(h)})$, observing that
\begin{eqnarray*}
&&
\int_{0}^{1} \frac
{e^{-(u_{i}^{2}-2r_{i,j}^{(h)}u_{i}u_{j}+u_{j}^{2})/2(1-(r_{i,j}^{(h)})^{2})}}{\sqrt
{1-(r_{i,j}^{(h)})^{2}}} \,dh \\[-2pt]
&&\qquad \leq \int_{0}^{1} \frac{1}{\sqrt
{1-(r_{i,j}^{(h)})^{2}}}
e^{-(u_{i}^{2}+u_{j}^{2})/2(1+|r_{i,j}^{(h)}|)} \,dh \\[-2pt]
&&\qquad = \frac{1}{r_{i,j}^{(1)} - r_{i,j}^{(0)}} \int
_{r_{i,j}^{(0)}}^{r_{i,j}^{(1)}} \frac{1}{\sqrt{1-t^{2}}}
e^{-(u_{i}^{2}+u_{j}^{2})/2(1+|t|)} \,dt.
\end{eqnarray*}
For bound (ii), overestimate the exponential by
$e^{-(u_{i}^{2}+u_{j}^{2})/2(1+\max\{|r_{i,j}^{(1)}|,|r_{i,j}^{(0)}|\}
)}$, and then evaluate the integral over $t$. Alternatively, by making
a substitution $x=\sqrt{(1-t)/(1+t)}$ we find that for any $0 \leq a
\leq b < 1$, and any $K \geq0$,
\begin{eqnarray*}
\int_{a}^{b} \frac{1}{\sqrt{1-t^2}} e^{-K/(1+t)} \,dt & = & 2e^{-K/2}
\int_{\sqrt{(1-b)/(1+b)}}^{\sqrt{(1-a)/(1+a)}} \frac{1}{1+x^{2}}
e^{-Kx^{2}/2} \,dx \\[-2pt]
& \leq& \frac{(1+b)^{3/2}}{\sqrt{1-b} K} e^{-K/2} \int_{\sqrt
{(1-b)/(1+b)}}^{\sqrt{(1-a)/(1+a)}} Kxe^{-Kx^{2}/2} \,dx.
\end{eqnarray*}
Since this integral is at most $e^{-K(1-b)/2(1+b)}$, the third bound
follows directly.\vspace*{1pt}

As Leadbetter, Lindgren and Rootz\'{e}n~\cite{llr} point out, the
assumption that $\tilde{X}$ and $\tilde{W}$ are nonsingular is not
necessary for the above bounds, as one may pass to that case by making
arbitrarily small changes to the entries of the covariance matrices,
and the first bound (from which we derived the others) is a continuous
function of those entries.

Typically, one would apply Comparison Inequality~\ref{comp1} by
observing that
the covariance matrix of $\tilde{X}$ ``looks rather like'' the
covariance matrix of a well understood multivariate normal
distribution, for example, that it looks like the identity matrix (see
the paper of Li and Shao~\cite{lishao} for some examples). If the
entries of the covariance matrices are sufficiently close together, or
if one can afford to choose the entries of $\tilde{u}$ very large,
then Comparison Inequality~\ref{comp1} can supply strong information.\vadjust{\goodbreak}

We finish with a well-known qualitative consequence of Comparison
Inequality~\ref{comp1}.
%
%
\begin{compineq}\label{comp2}
Let $\tilde{X}=(X_{1},\ldots,X_{n})$, $\tilde{W}=(W_{1},\ldots,W_{n})$ be
centralized and normalized $n$-variate normal vectors, with covariance
matrices $\var(\tilde{X})= (r_{i,j}^{(1)})$ and $\var(\tilde{W})=
(r_{i,j}^{(0)})$, respectively. Let $\tilde{u} \in\R^{n}$. If
$r_{i,j}^{(1)} \leq r_{i,j}^{(0)}$ for each $1 \leq i,j \leq n$, then
\[
\p(\tilde{X} \leq\tilde{u}) \leq\p(\tilde{W} \leq\tilde{u}).
\]
\end{compineq}

The special case of this result where $\tilde{u}=(u,u,\ldots,u)$, for
some $u \in\R$, is usually referred to as Slepian's lemma.

\subsection{Reversal of roles}\label{sec3.2}
As promised in the \hyperref[sec1]{Introduction}, we now give a very brief description
of the reversal of roles argument that originally served in place of
Proposition~\ref{pr1}. For the applications in this paper, Proposition
\ref{pr1}
entirely supersedes such an argument, but it is possible that it may be
useful in other contexts.

We aim to give an estimate for
\[
\p\bigl(Z(t_{m})>u, Z(t_{j}) \leq u\ \forall j \leq m-1\bigr)
\]
under the conditions of Proposition~\ref{pr1}. Our idea is to apply the
methodology of Exact Formula~\ref{ex1}, but viewing the sum of
integrals that
arises as a main term for subsequent analysis and the subtracted
probability as an error term. Thus we do not choose $\tilde{W}$ to
have a standard distribution, but so that this subtracted probability
is zero.

More concretely, we let $A_{1},\ldots,A_{m}$ be a collection of $N(0,1)$
random variables, all independent of one another and of the $Z(t_{i})$.
Let $\varepsilon> 0$, and define
\begin{eqnarray*}
X_{i} &=& W_{i}:= \frac{Z(t_{i})+\varepsilon A_{i}}{\sqrt{1+\varepsilon
^{2}}},\qquad 1 \leq i \leq m-1;
\\
X_{m}:\!&=& \frac{Z(t_{m}) +\varepsilon
A_{m}}{\sqrt{1+\varepsilon^{2}}};\qquad
W_{m}:= \frac{Z(t_{m-1}) + \varepsilon A_{m}}{\sqrt{1+\varepsilon^{2}}}.
\end{eqnarray*}
Precisely analogously to Exact Formula~\ref{ex1}, and adopting the same
notation $r_{i,j}^{(h)}$ as there, we find that
\begin{eqnarray*}
&&\p(X_{1},\ldots,X_{m-1} \leq u, X_{m} > u) - \p(W_{1},\ldots,W_{m-1}
\leq
u, W_{m} > u)
\\
&&\qquad= - \sum_{1 \leq i \leq m-1} \bigl(r_{i,m}^{(1)} - r_{i,m}^{(0)}\bigr)
\int
_{0}^{1} \phi\bigl(u,u; r_{i,m}^{(h)}\bigr)\\
&&\hspace*{149.2pt}{}\times \p\bigl(\tilde{Z}^{(h)} \leq\tilde{u}
| Z^{(h)}_{i}=u, Z^{(h)}_{m}=u\bigr) \,dh
\\
&&\qquad= - \sum_{1 \leq i \leq m-1} \frac
{r_{i,m}-r_{i,m-1}}{1+\varepsilon
^{2}} \int_{0}^{1} \phi\biggl(u,u; \frac
{hr_{i,m}+(1-h)r_{i,m-1}}{1+\varepsilon^{2}}\biggr) \\
&&\hspace*{156pt}{}\times P(i,h,\varepsilon) \,dh,
\end{eqnarray*}
say. We need the $\varepsilon$ perturbations here to ensure that we work
with nonsingular multivariate normal distributions. However, at the end
of the argument we can let $\varepsilon\rightarrow0$, whereby we will
have compared $\p(Z(t_{m})>u, Z(t_{j}) \leq u\ \forall j \leq m-1)$ with
\[
\p\bigl(Z(t_{1}),\ldots,Z(t_{m-1}) \leq u, Z(t_{m-1}) > u\bigr)=0.
\]

It is less straightforward to analyze $P(i,h,\varepsilon)$ for $1 \leq i
\leq m-2$ than to analyze $P(m-1,h,\varepsilon)$, and to give lower
bounds one can replace those probabilities by $\mathbf
{1}_{r_{i,m}>r_{i,m-1}}$. In our examples, these other terms give a
lower order contribution, but this need not always be so. However, to
analyze $P(m-1,h,\varepsilon)$ one can note (as did Li and Shao \cite
{lishao}) that for any $1 \leq i \leq m-1$, the collection of random variables
\begin{eqnarray}
Y_{j}^{(h)}:\!&=& Z_{j}^{(h)}-\biggl(\frac{r_{j,i}^{(h)} -
r_{i,m}^{(h)}r_{j,m}^{(h)}}{1-(r_{i,m}^{(h)})^{2}} \biggr)Z_{i}^{(h)} -
\biggl(\frac{r_{j,m}^{(h)} -
r_{i,m}^{(h)}r_{j,i}^{(h)}}{1-(r_{i,m}^{(h)})^{2}} \biggr)Z_{m}^{(h)}
\nonumber\\
&=& Z_{j}^{(1)} - \biggl(\frac{r_{j,i}^{(1)} -
r_{i,m}^{(h)}r_{j,m}^{(h)}}{1-(r_{i,m}^{(h)})^{2}} \biggr)Z_{i}^{(1)} -
\biggl(\frac{r_{j,m}^{(h)} -
r_{i,m}^{(h)}r_{j,i}^{(1)}}{1-(r_{i,m}^{(h)})^{2}}
\biggr)Z_{m}^{(h)},\nonumber\\
&&\eqntext{1 \leq j \leq m-1, j \neq i,}
\end{eqnarray}
is independent of $\{Z_{i}^{(h)},Z_{m}^{(h)}\}$. In our examples this
leads, after some slightly fiddly manipulations, to a probability
estimate much like Proposition~\ref{pr1}. [Since, in our examples,
$Z(t_{m-1})$ and
$Z(t_{m})$ are always very highly correlated, and so
$P(m-1, h,\varepsilon
)$ is essentially the same as the simple conditional probability in the
proof of Proposition~\ref{pr1}.]

\section{\texorpdfstring{Proof of Proposition \protect\ref{pr2}}{Proof of Proposition 2}}\label{sec4}
In view of Comparison Inequality~\ref{comp2}, and assumption (ii) in the
statement of Proposition~\ref{pr2}, we may proceed on the assumption
that for
$1 \leq j,k \leq m-1$ and $j \neq k$, $\E V_{j}V_{k}$ is equal to
\[
\frac{c_{\min\{j,k\}} d_{\max\{j,k\}}}{\sqrt{(1 - r_{j,m}^{2})(1 -
r_{k,m}^{2})}}.
\]
The key to the proof is the explicit construction of such random
variables from a collection of independent normal random variables.

Let $Y_{1},\ldots,Y_{n},Z_{1},\ldots,Z_{n}$ be independent standard normal
random variables, and for $1 \leq i \leq n$ let $\alpha_{i}, \beta
_{i}$ be real numbers satisfying
\[
\beta_{i}^{2} \sum_{j \leq i} \alpha_{j}^{2} < 1.
\]
Then the random variables
\[
X_{i}:= \beta_{i} \sum_{j \leq i} \alpha_{j} Y_{j} + \sqrt{1-\beta
_{i}^{2} \sum_{j \leq i} \alpha_{j}^{2}} Z_{i}
\]
are again jointly multivariate normal, have zero means and unit
variances and satisfy
\[
\E X_{i} X_{j} = \beta_{i} \beta_{j} \sum_{k \leq\min\{i,j\}}
\alpha_{k}^{2},\qquad i \neq j.
\]
We also note that if $u_{1},\ldots,u_{n}$ are any real numbers, if
$\beta
_{i} > 0\ \forall1 \leq i \leq n$ and if $\delta\in\R$, then
\begin{eqnarray*}
&&\p(X_{i} \leq u_{i}\ \forall1 \leq i \leq n) \\
&&\qquad = \p\biggl(Z_{i} \leq
\frac{u_{i}-\beta_{i} \sum_{j \leq i} \alpha_{j} Y_{j}}{\sqrt
{1-\beta_{i}^{2} \sum_{j \leq i} \alpha_{j}^{2}}}\ \forall1 \leq i
\leq n \biggr) \\
&&\qquad \geq \p\biggl(\sum_{j \leq i} \alpha_{j} Y_{j} \leq\frac{\delta
u_{i}}{\beta_{i}}\ \forall1 \leq i \leq n\biggr) \prod_{i=1}^{n} \Phi
\biggl(\frac{u_{i}(1-\delta)}{\sqrt{1-\beta_{i}^{2} \sum_{j \leq i}
\alpha_{j}^{2}}} \biggr).
\end{eqnarray*}

We now set $n=m-1$, and define real numbers $\alpha_{i},\beta_{i}$ by
\[
\beta_{i}:= \frac{d_{i}}{\sqrt{1 - r_{i,m}^{2}}},\qquad \sum_{j \leq i}
\alpha_{j}^{2}:= \frac{c_{i}}{d_{i}},\qquad 1 \leq i \leq m-1.
\]
The conditions on $c_{i},d_{i}$ in Proposition~\ref{pr2} ensure that we
can define $\alpha_{i},\beta_{i}$ in this way, and that they satisfy
the various hypotheses above. The reader may also check that the
$X_{i}$ have the correlation structure that we wanted, and that the
product term in the previous paragraph is as in
Proposition~\ref{pr2} [when $u_{i}$ is taken as
$(u-r_{i,m}(u+h))/\sqrt{1-r_{i,m}^{2}}$]. It remains to give a suitable
lower bound for $\p(\sum_{j \leq i} \alpha_{j} Y_{j} \leq \frac{\delta
u_{i}}{\beta_{i}}\ \forall1 \leq i \leq m-1)$.\vspace*{1pt}

It should not come as a surprise that the behavior of partial sums of
independent normal random variables is rather well understood. For
example, writing $\{W_{t}\}_{t \geq0}$ for the standard Brownian
motion (see, e.g., Lifshits~\cite{lif}, Chapter 5, for much discussion
of this process), one has the following neat result, which we quote
from Grimmett and Stirzaker~\cite{gs}, Chapter 13.4: if $t \geq0$, then
\[
\max_{0 \leq s \leq t} W_{s} \stackrel{d}{=} |W_{t}| \stackrel{d}{=}
|N(0,t)|.
\]
This is useful to us because $(\sum_{j \leq i} \alpha_{j} Y_{j} )_{1
\leq i \leq m-1} \stackrel{d}{=} (W_{\sum_{j \leq i} \alpha
_{j}^{2}})_{1 \leq i \leq m-1}$, so that
\[
\p\biggl(\sum_{j \leq i} \alpha_{j} Y_{j} \leq\frac{\delta u_{i}}{\beta
_{i}}\ \forall1 \leq i \leq m-1\biggr) \geq\Phi(B)-\Phi(-B),
\]
where
\[
B = \frac{\delta}{\sqrt{\sum_{j \leq m-1} \alpha_{j}^{2}}} \min
_{1 \leq i \leq m-1} \frac{u_{i}}{\beta_{i}} = \delta\sqrt{\frac
{d_{m-1}}{c_{m-1}}} \min_{1 \leq i \leq m-1} \frac{u-r_{i,m}(u+h)}{d_{i}}
\]
as claimed in Proposition~\ref{pr2}.\qed\vspace*{8pt}
%

The proof just given divided naturally into two parts: first we
constructed the $X_{j}$ to explicitly model the $V_{j}$, allowing us to
extract some of their dependence in the manageable form of the $Y_{j}$;
and then we analyzed the $Y_{j}$ using a result about Brownian motion.
Both of these steps could conceivably be improved, potentially leading
to a better lower bound for $P(m,h)$.

In the analysis of the $Y_{j}$, we used a fact about the probability
that a Brownian motion remains below a constant level for a period of
``time'' $t$. We could have used results about the probability that it
remains below, for example, a sloping line, allowing some flexibility
in the upper bounds that we ask for. However, in our applications these
probabilities are never particularly small, and the author doubts that
a more complicated approach would be advantageous in many situations.

It appears to the author that the modeling part of the argument is
weaker. Thus, in our examples, our lower bound $c_{\min\{j,k\}}d_{\max
\{j,k\}}$ for $r_{j,k}-r_{j,m}r_{k,m}$ is not very tight when $j$ and
$k$ are close together. An alternative way to think about this is to
note that we can replace the independent $Z_{j}$ in our construction by
any standard normal $A_{j}$ with
\[
\E A_{j}A_{k} \leq\frac{r_{j,k}-r_{j,m}r_{k,m}-c_{\min\{j,k\}
}d_{\max\{j,k\}}}{\sqrt
{(1-r_{j,m}^{2}-c_{j}d_{j})(1-r_{k,m}^{2}-c_{k}d_{k})}}.
\]
The correlation bound here looks complicated, but this may be somewhat
illusory; for example, if we were able to make the choices
$c_{j}=r_{j,m}$, $d_{j}=1-r_{j,m}$, as for certain stationary
processes, we would want
\[
\E A_{j}A_{k} \leq\frac{r_{j,k}-r_{\min\{j,k\},m}}{\sqrt
{(1-r_{j,m})(1-r_{k,m})}}.
\]
These quantities are not likely to be easier to work with than the
correlations $r_{j,k}$ of our original random variables. However, to
prove Proposition~\ref{pr2} we need upper bounds for upper tail probabilities
(which then lower bound the probability that none of the $A_{j}$ are
too big), and these may be easier to come by than lower bounds, for
example by using Rice's formula as part of a first moment argument.
Another approach to improving Proposition~\ref{pr2} along these lines
is worked
out in Section~\ref{sec7}.

\section{Application to estimating Pickands' constants}\label{sec5}
Suppose that $t_{1} < t_{2} <\cdots< t_{M}$ is a set of equally spaced
real numbers. Suppose, moreover, that $\{Z(t_{i})\}_{1 \leq i \leq M}$
is a mean zero, variance one, stationary Gaussian process with
decreasing covariance function $r(t)$, $t \geq0$. If $a > 0$, then
(the proof of) Proposition~\ref{pr1} implies that
\begin{eqnarray*}
\p\Bigl(\max_{1 \leq i \leq M} Z(t_{i}) > u\Bigr) & \geq& M \p\bigl(Z(t_{M}) > u,
Z(t_{j}) \leq u\ \forall j < M\bigr) \\
& \geq& \frac{M e^{-u^{2}/2}}{\sqrt{2\pi} u} \cdot a
e^{-a-a^{2}/2u^{2}} \inf_{0 \leq h \leq a/u} P(M,h).
\end{eqnarray*}

In a paper from 1996, Shao~\cite{shao} considers a mean zero, variance
one stationary Gaussian process indexed by the half-line $[0,\infty)$,
with covariance function
\[
r(t)=\tfrac{1}{2}\bigl(e^{\alpha t/2} + e^{-\alpha t/2}
-(e^{t/2}-e^{-t/2})^{\alpha}\bigr),\qquad t \geq0.
\]
Such a process exists for each fixed $0 < \alpha< 2$. As $t
\rightarrow0$, we see (as did Shao~\cite{shao}) that
$r(t)=1-t^{\alpha}/2+O(t^{2})$. We also note that, for $t > 0$,
\begin{eqnarray*}
r'(t) & = & \frac{\alpha}{4}\bigl(e^{\alpha t/2}-e^{-\alpha t/2}
-(e^{t/2}+e^{-t/2})(e^{t/2}-e^{-t/2})^{\alpha-1}\bigr) \\
& = & \frac{\alpha}{4}\bigl(e^{\alpha t/2}-e^{-\alpha t/2} -e^{\alpha
t/2}(1+e^{-t})(1-e^{-t})^{\alpha-1}\bigr) \\
& \leq& \frac{\alpha}{4}\bigl(e^{\alpha t/2}-e^{-\alpha t/2} -e^{\alpha
t/2}(1-e^{-2t})\bigr)
< 0.
\end{eqnarray*}

From now on it will be convenient to employ Vinogradov's notation $\gg
$, meaning ``greater than, up to a multiplicative constant.'' Thus
$p(\alpha) \gg q(\alpha)$ means the same as $q(\alpha) = O(p(\alpha
))$. In proving Corollary~\ref{co1}, we shall assume that $\alpha$ is smaller
than a certain positive constant less than 1; a suitable explicit value
could be extracted from our calculations if desired. There is no loss
in this because $H_{\alpha} \gg1 \gg\sqrt{\alpha} (e \alpha
/2)^{1/\alpha}$ for $\alpha$ larger than such a constant. To prove
the corollary, we will study Shao's stationary process at the sample
points $t_{i}=i/M$. Simply choosing $a=1$ in the above discussion, and
comparing with Pickands' theorem in the \hyperref[sec1]{Introduction}, we see
\[
H_{\alpha} \gg2^{1/\alpha} \lim_{u \rightarrow\infty} \Bigl(M
u^{-2/\alpha} \inf_{0 \leq h \leq1/u} P(M,h)\Bigr),
\]
and we will investigate the largest value of $M$, depending on $u$ and
$\alpha$, for which we can show that $\inf_{0 \leq h \leq1/u} P(M,h)
\gg1$. Note that a large value of $M$ corresponds to a close packing
of sample points in the interval $[0,1]$. The reader should also note
that there is nothing intrinsically asymptotic about most of our
calculations, although we are interested in letting $u \rightarrow
\infty$ to compare with Pickands' theorem.

We want to apply Proposition~\ref{pr2}, and can do so with the natural choices
\[
c_{j}=r\bigl((M-j)/M\bigr),\qquad d_{j}=1-r\bigl((M-j)/M\bigr),
\qquad 1 \leq j \leq M-1,
\]
since $r(t)$ is decreasing and positive. Thus $P(M,h)$ is at least
\begin{eqnarray*}
&&\bigl(\Phi(B)-\Phi(-B)\bigr) \prod_{j=1}^{M-1} \Phi\biggl(\biggl(1+O\biggl(\frac
{1}{u^{2}(1-r(j/M))}\biggr)\biggr)u(1-\delta)\sqrt{1-r(j/M)} \biggr)
\\
&&\qquad= \bigl(\Phi(B)-\Phi(-B)\bigr)\\
&&\qquad\quad{}\times
\prod_{j=1}^{M-1} \Phi\biggl(\biggl(1+O\biggl(\frac{M^{\alpha
}}{u^{2} j^{\alpha}}\biggr)\biggr)(1-\delta)\sqrt{u^{2}\bigl(j^{\alpha}/2M^{\alpha}
+ O(j^{2}/M^{2})\bigr)} \biggr),
\end{eqnarray*}
where $B=B(\delta)$ is as in Proposition~\ref{pr2}, and $\delta$ will be
chosen later in terms of~$\alpha$. Together with the known and
conjectured bounds for Pickands' constants, this suggests taking
$M=[(bu^{2}\alpha/2)^{1/\alpha}]$, where now we investigate how large
$b$ may be chosen. For definiteness in our calculations, we declare
that we shall certainly have $1 \leq b \leq10$ (and of course our
conclusion will be that taking $b$ as $e/2$ is permissible).

First we note that the part of the product over $j > M^{1/4}$ is
$1+o(1)$ as $u \rightarrow\infty$. For, since $r(t)$ is decreasing
and $\alpha,\delta$ are small, each of those terms is at least
\begin{eqnarray*}
&&
\Phi\bigl(\bigl(1+O(\alpha)\bigr)u(1-\delta)\sqrt{1-r(M^{-3/4})}\bigr)\\
&&\qquad \geq\Phi
\bigl((1/2)\sqrt{u^{2}\bigl(M^{-3\alpha/4}/2+O(M^{-3/2})\bigr)}\bigr).
\end{eqnarray*}
If $u$, and therefore the argument of $\Phi$, is large enough, this is
\[
\geq1 - e^{-(1/8)u^{2}(M^{-3\alpha/4}/2 + O(M^{-3/2}))} \geq1 -
e^{-\sqrt{u}},
\]
and clearly $(1-e^{-\sqrt{u}})^{M}$ is $1+o(1)$ as $u \rightarrow
\infty$ with $\alpha$ fixed.

When $j \leq M^{1/4}$, provided that $u$ is large enough in terms of
$\alpha\leq1$ we see
\[
u^{2}j^{2}/M^{2} \leq u^{2} M^{-3/2} = O(u^{-1}\alpha^{-3/2})
\]
and
\[
M^{\alpha}/u^{2}j^{\alpha} = O(\alpha/j^{\alpha})
\]
and so the terms in the product are
\begin{eqnarray*}
&&
\Phi\bigl(\bigl(1+O(\alpha/j^{\alpha})\bigr)(1-\delta)\sqrt{j^{\alpha}/(b\alpha
) + O(1/u\alpha^{3/2})} \bigr)\\
&&\qquad = \Phi\bigl(\bigl(1+O(\alpha/j^{\alpha})\bigr)(1-\delta
)\sqrt{j^{\alpha}/(b\alpha)} \bigr).
\end{eqnarray*}
Thus, since $\Phi(x) \geq1 - x^{-1}e^{-x^{2}/2} \geq
e^{-2x^{-1}e^{-x^{2}/2}} $ for $x \geq2$, the part of the product over
$j \leq M^{1/4}$ is at least $e^{-f(b,\delta,\alpha,u)}$, where
\[
f(b,\delta,\alpha,u) = O\biggl(\sum_{j \leq M^{1/4}} e^{-(1-\delta
)^{2}j^{\alpha}/2b\alpha} \sqrt{b\alpha}/(1-\delta)j^{\alpha/2}\biggr).
\]
(Since we assume that $\alpha$ and $\delta$ are small, the arguments
of $\Phi$ are all at least~$2$.) Now
\begin{eqnarray*}
\sum_{j \leq M^{1/4}} e^{-(1-\delta)^{2}j^{\alpha}/2b\alpha} & \leq
& \int_{0}^{M^{1/4}} e^{-(1-\delta)^{2}t^{\alpha}/2b\alpha} \,dt
\\
& = & \frac{2b}{(1-\delta)^{2}} \int_{0}^{(1-\delta)^{2} M^{\alpha
/4}/2b\alpha} \biggl( \frac{2b \alpha y}{(1-\delta)^{2}} \biggr)^{1/\alpha- 1}
e^{-y} \,dy \\
& \leq& \biggl(\frac{2b}{(1-\delta)^{2}}\biggr)^{1/\alpha} \alpha^{1/\alpha-
1} \Gamma(1/\alpha).
\end{eqnarray*}
By Stirling's formula, the right-hand side is asymptotic to
\[
\sqrt
{2\pi/\alpha} \bigl(2b/e(1-\delta)^{2}\bigr)^{1/\alpha}
\]
as $\alpha \rightarrow0$, so is at most $4 \alpha^{-1/2}
(2b/e(1-\delta)^{2})^{1/\alpha}$, say, when $\alpha$ is small.

Finally, observe that
\[
B(\delta) = \delta\sqrt{\bigl(1-r(1/M)\bigr)/r(1/M)} u
\bigl(1+O\bigl(1/u^{2}\bigl(1-r(1/M)\bigr)\bigr)\bigr) \gg\delta/\sqrt{\alpha},
\]
provided that $u$ is large enough in terms of $\alpha$. If we make the
choice $\delta= \alpha$, then $b$ can be chosen as large as $e/2$
while still ensuring that $f(b,\delta,\alpha,u)=O(1)$. Corollary
\ref{co1}
follows from making these choices.\qed

\section{Application to a number-theoretic process}\label{sec6}

\subsection{Preliminary calculations}\label{sec6.1}
Before we can apply Propositions~\ref{pr1} and~\ref{pr2} to our second
example, we must
reduce to studying a finite set of sample points $t$ and determine the
covariance structure of the corresponding random variables. As might be
expected, variants of some of these calculations appear in Hal\'{a}sz's
paper~\cite{hal}, but we must be more precise in several places.

It is useful initially to ignore the contribution from ``very small''
primes to our random sums. Let $y$ be a parameter, later to be chosen
as a suitable function of $x$. It is immediate that if $s,t \in\R$, then
\begin{eqnarray*}
&& \E\biggl(\sum_{y \leq p \leq x} g_{p} \frac{\cos(t \log
p)}{p^{1/2+1/\log x}} \cdot\sum_{y \leq p \leq x} g_{p} \frac{\cos
(s \log p)}{p^{1/2+1/\log x}}\biggr) \\
&&\qquad= \sum_{y \leq p \leq x} \frac{\cos(t\log p) \cos(s\log
p)}{p^{1+2/\log x}} \\
&&\qquad= \frac{1}{2} \sum_{y \leq p \leq x} \frac{\cos((t+s)\log p) +
\cos((t-s)\log p)}{p^{1+2/\log x}}.
\end{eqnarray*}
For $t \in\R$, we let $Z_{y}(t)$ denote the normalized random variable
\begin{eqnarray*}
&&\frac{\sum_{y \leq p \leq x} g_{p}\cos(t \log p)/p^{1/2+1/\log
x}}{\sqrt{\sum_{y \leq p \leq x} \cos^{2}(t \log p)/p^{1+2/\log x}}}\\
&&\qquad
= \frac{\sum_{y \leq p \leq x} g_{p}\cos(t \log p)/p^{1/2+1/\log
x}}{\sqrt{(\sum_{y \leq p \leq x} 1/p^{1+2/\log x} +
\sum_{y \leq p \leq x} \cos(2t \log p)/p^{1+2/\log x})/2}}.
\end{eqnarray*}

By a strong form of the prime number theorem (see, e.g., Montgomery and
Vaughan~\cite{mv}, Chapter 6) we have
\begin{eqnarray*}
\pi(z):= \#\{p \leq z\dvtx p \mbox{ is prime}\}
= \int_{2}^{z} \frac
{du}{\log u} + O\bigl(ze^{-d\sqrt{\log z}}\bigr),\qquad z \geq2,
\end{eqnarray*}
where $d > 0$ is a certain constant. Then if $z \leq x$,
\begin{eqnarray*}
\sum_{p \leq z} \frac{1}{p^{1+2/\log x}} &=& \int_{2}^{z} \frac
{1}{u^{1+2/\log x}} \,d\pi(u) \\
& = & \int_{2}^{z} \frac{u^{-2/\log
x}}{u\log u}\,du + c(x) + O\bigl(e^{-d\sqrt{\log z}}\bigr) \\
& = & \log\log z + O(1),
\end{eqnarray*}
where $c(x)$ depends on $x$ only. Moreover, if $\alpha\neq0$,
\begin{eqnarray*}
\sum_{y \leq p \leq x} \frac{\cos(\alpha\log p)}{p^{1+2/\log x}}
& = & \int_{y}^{x} \frac{\cos(\alpha\log u) u^{-2/\log x}}{u\log u}
\,du
+ O\bigl((1+|\alpha|)e^{-d\sqrt{\log y}}\bigr) \\
& = & \int_{\log y}^{\log x} \frac{\cos(\alpha u)}{u} \,du + \int
_{\log y}^{\log x} \frac{\cos(\alpha u)}{u}(e^{-2u/\log x} - 1) \,du
\\
&&{} + O\bigl((1+|\alpha|)e^{-d\sqrt{\log y}}\bigr) \\
& = & \int_{\alpha\log y}^{\alpha\log x} \frac{\cos u}{u} \,du
+
O\biggl(\frac{1}{\alpha\log x}\biggr)
+ O\bigl((1+|\alpha|)e^{-d\sqrt{\log y}}\bigr),
\end{eqnarray*}
where the third equality follows using integration by parts since
$\frac{d}{du}((e^{-2u/\log x}-1)/u) = O(1/\log^{2}x)$ for $\log y
\leq u \leq\log x$. We deduce that if $s,t \geq1$, and $s \neq t$, then
\begin{eqnarray*}
\E Z_{y}(t)Z_{y}(s)
&=& \biggl(\int_{|t-s|\log y}^{|t-s|\log x} \frac
{\cos u}{u}\,du + O\biggl(\frac{1}{(t+s)\log y}\biggr)\\
&&\hspace*{4.7pt}{} + O\biggl(\frac{1}{|t-s|\log x}\biggr) +
O\bigl((t+s)e^{-d\sqrt{\log y}}\bigr)\biggr)\\
&&{}\times\biggl({\int_{y}^{x} \frac{du}{u^{1+2/\log x}
\log u} + O\biggl(\frac{1}{\log y}\biggr) + O\bigl((t+s)e^{-d\sqrt{\log y}}\bigr)}\biggr)^{-1}.
\end{eqnarray*}

We now set out the specific situation to which our Gaussian process
results will be applied. Let $E \geq1$ be a further parameter (to be
chosen later as a function of $x$) and for $n \in\N\cup\{0\}$ and $M
\leq(\log x)/E$ introduce the sets
\[
\mathcal{T}_{n}=\mathcal{T}_{n,x,E,M}:= \{2n+1+iE/\log x\dvtx1 \leq i
\leq M\} \subseteq[2n+1,2n+2].
\]
We seek lower bound information on $\max_{0 \leq n \leq B} \sup_{t
\in\mathcal{T}_{n}} Z_{y}(t)$, for certain $B$.

At this point the reader may be rather appalled by the number of
parameters around, so we hasten to point out that most of these will
``select themselves'' in due course and can essentially be ignored. The
sets $\mathcal{T}_{n}$ are sufficiently separated that the behavior of
$Z_{y}(t)$ on different sets is roughly independent (see
Section~\ref{sec6.3}). Moreover, up to error terms the correlation $\E
Z_{y}(t)Z_{y}(s)$ depends on $s,t$ through $|t-s|$ only (i.e., our
process is approximately stationary). Thus
we focus
on understanding $\sup_{t \in \mathcal{T}_{0}} Z_{y}(t)$, and defer
thinking about larger values of $n$ until we put our results together
in Section~\ref{sec6.3}.

The parameter $E$ controls the spacing of sample points within their
blocks $\mathcal{T}_{n}$, and in Section~\ref{sec6.2} it will be chosen as
small as possible such that we obtain good probability lower bounds
from Proposition~\ref{pr2}. We declare for now that we shall certainly
have $E
\leq e^{\sqrt{\log\log x}}$, say. We would like to take $M$ as large
as possible, but to simplify our calculations we choose $M=[\log x/KE
\log y]$, where $K$ is an absolute constant that forces $\E
Z_{y}(t)Z_{y}(s) \geq1/\log\log x$, say, for $t,s \in\mathcal
{T}_{0}$ (see below). The parameter $y$ is present to get rid of
``beginning of series'' effects, in particular ensuring that we have
enough independence of $Z_{y}(t)$ over different blocks $\mathcal
{T}_{n}$. It will be selected in Section~\ref{sec6.3}, but we declare
for now
that we shall certainly have $\log x \leq y \leq e^{(\log\log x)^{100}}$.

In the above set-up, if $s,t \in\mathcal{T}_{0}$ are distinct, then
\begin{eqnarray*}
\E Z_{y}(t)Z_{y}(s) & = & \frac{\int_{|t-s|\log y}^{\log x}
({\cos u}/{u})\,du + O({1}/({|t-s|\log x}))}{\int_{y}^{x}
{du}/({u^{1+2/\log x} \log u})}\\
&&{} + O\biggl(\frac{1}{\log y \log\log x} \biggr)
\\
& = & \frac{\int_{|t-s|\log y}^{1} ({\cos u}/{u})\,du }{\int
_{y}^{x} {du}/({u^{1+2/\log x} \log u})} + O\biggl(\frac{1}{\log\log x}
\biggr) \\
& = & \frac{\log(1/|t-s|\log y)}{\log\log x - \log\log y} + O\biggl(\frac
{1}{\log\log x} \biggr).
\end{eqnarray*}

\subsection{\texorpdfstring{Implementation of Propositions \protect\ref{pr1} and \protect\ref{pr2}}
{Implementation of Propositions 1 and 2}}\label{sec6.2}
We order the points of $\mathcal{T}_{0}$ in the obvious and natural
way, writing $t_{i} = 1+iE/\log x$, $1 \leq i \leq M$. We aim to show
that the maximum of our original random sum is about $\log\log x$, and
the standard deviations that we normalized by are about $\sqrt{(\log
\log x - \log\log y)/2}$, so we take $u=\sqrt{2(\log\log x -\log
\log y)}$. Then, recalling our notation $r_{m-1,m} = \E
Z_{y}(t_{m-1})Z_{y}(t_{m})$,
\[
u(1-r_{m-1,m})=\Theta\bigl(\log E/\sqrt{\log\log x - \log\log y}\bigr)=\Theta
(\log E/ u),
\]
so we can safely make the canonical choice $H=1/u$ in Proposition
\ref{pr1}.

We now seek to apply Proposition~\ref{pr2} to give a lower bound for $P(m,h)$,
where $1 \leq m \leq M$ and $h \leq H$. Let $j < k \leq m-1$. If $|j-k|
\leq\log^{1/3}x$ then
\begin{eqnarray*}
r_{j,k} &=& 1 - \frac{\log(|j-k|E)}{\log\log x - \log\log y} +
O\biggl(\frac{1}{\log\log x}\biggr) \\
&\geq& \max\{1/2,r_{j,m}\} + O\biggl(\frac
{1}{\log\log x}\biggr) \\
& \geq& r_{j,m} + O\biggl(\frac{r_{j,m}}{\log\log x}\biggr).
\end{eqnarray*}
In fact this is also true when $|j-k| > \log^{1/3}x$. Since $\int
_{\alpha\log y}^{\log x} (\cos u/u) \,du$ is a decreasing function of $0
< \alpha\leq1/\log y$, we have
\begin{eqnarray*}
r_{j,k} & = & \frac{\int_{|j-k|E\log y /\log x}^{\log x} ({\cos
u}/{u})\,du + O({1}/({|j-k|E}))}{\int_{y}^{x} {du}/({u^{1+2/\log x}
\log u})} \\
&& {}+ O\biggl(\frac{1}{\log y \log\log x} \biggr) \\
& \geq& r_{j,m} + O\biggl(\frac{1}{\log y \log\log x} \biggr)
\end{eqnarray*}
and we always have $r_{j,m} \geq1/\log\log x \geq1/\log y$ because
$|m-j| \leq M$. This means that $r_{j,k} - r_{j,m}r_{k,m} \geq
r_{j,m}(1-r_{k,m}+O(1/\log\log x))$, so it is legitimate to choose
\begin{eqnarray*}
c_{j}&=&1-\frac{\log((m-j)E) + O(1)}{\log\log x - \log\log y},\\
d_{j}&=&\frac{\log((m-j)E) + O(1)}{\log\log x - \log\log y},\qquad 1 \leq j
\leq m-1,
\end{eqnarray*}
in Proposition~\ref{pr2}. Setting $\delta=1/\log\log x$ in the proposition,
to match the size of our other ``big Oh'' terms, we discover that
\[
B(1/\log\log x)= \Theta\biggl(\frac{u \sqrt{\log E}}{(\log\log x)^{3/2}}
\biggr) = \Theta\biggl(\frac{\sqrt{\log E}}{\log\log x} \biggr).
\]
It follows from all of this that, for $1 \leq m \leq M$ and $h \leq H$,
\begin{eqnarray*}
P(m,h) & \gg& \frac{\sqrt{\log E}}{\log\log x} \prod_{j=1}^{m-1}
\Phi\biggl( \biggl(1+O\biggl(\frac{1}{\log((m-j)E)}\biggr)\biggr) \sqrt{2\log\bigl((m-j)E\bigr)} \biggr)
\\
& \gg& \frac{\sqrt{\log E}}{\log\log x} e^{-\Theta(\sum
_{j=1}^{m-1} {1}/({(m-j)E \sqrt{\log((m-j)E)}}))},
\end{eqnarray*}
provided always that $E$ is larger than an absolute constant. Making
the choice $E=\sqrt{\log\log x}$, the exponential becomes $\Theta
(1)$, and we find $P(m,h) \gg\sqrt{\log\log\log x}/\log\log x$.

Plugging this lower bound into Proposition~\ref{pr1}, it follows
immediately that
\begin{eqnarray*}
\p\Bigl(\sup_{t \in\mathcal{T}_{0}} Z_{y}(t) > \sqrt{2(\log\log x -
\log\log y)}\Bigr) & \gg& \frac{M \sqrt{\log\log\log x}
e^{-u^{2}/2}}{u \log\log x} \\
& \gg& \frac{\sqrt{\log\log\log x}}{(\log\log x)^{2}}.
\end{eqnarray*}

\subsection{Exploitation of the lower bound}\label{sec6.3}
The lower bound obtained at the end of Section~\ref{sec6.2} is useful raw
information about $\{Z_{y}(t)\}_{t \in\mathcal{T}_{0}}$. However, in
order to deduce results about the summatory function $M(x)$ of a random
multiplicative function, as described in the \hyperref[sec1]{Introduction}, we need to
be able to say that the supremum is large with probability close to 1.

To do this, our idea is to ``sample the supremum several times
independently.'' Since the probability that the supremum is large is
not too small, if we just sample a few times we will very likely obtain
a large value. Although we do not have lots of independent copies of $\{
Z_{y}(t)\}$, we can achieve something like this by considering $\{
Z_{y}(t)\}_{t \in\mathcal{T}_{n}}$ for different $n$. If $Be^{-d\sqrt
{\log y}} \leq\frac{1}{\log y}$, say, then for distinct $1 \leq s,t
\leq2B+2$ we have
\begin{eqnarray*}
\E Z_{y}(t)Z_{y}(s) &=& \frac{\int_{|t-s|\log y}^{\log x} ({\cos
u}/{u})\,du + O({1}/({|t-s|\log x}))}{\int_{y}^{x}
{du}/({u^{1+2/\log x} \log u})} \\
&&{}+ O\biggl(\frac{1}{\log y \log\log x} \biggr)
\end{eqnarray*}
as at the end of Section~\ref{sec6.1}. For such $s,t$ with $|s-t| \geq
1$, the
calculations in Section~\ref{sec6.1} supply a more precise result,
namely, that
\[
\E Z_{y}(t)Z_{y}(s) = O\biggl(\frac{1}{|t-s|\log y \log\log x} + \frac
{(t+s)e^{-d\sqrt{\log y}}}{\log\log x} \biggr).
\]
Thus, by the second bound in Comparison Inequality~\ref{comp1},
\begin{eqnarray*}
&& \biggl|\p\Bigl(\max_{0 \leq n \leq B} \sup_{t \in\mathcal{T}_{n}} Z_{y}(t)
\leq\sqrt{2(\log\log x - \log\log y)}\Bigr) \\
&&\quad\hspace*{0pt}{} - \prod_{0 \leq n \leq B} \p\Bigl(\sup_{t \in\mathcal{T}_{n}}
Z_{y}(t) \leq\sqrt{2(\log\log x - \log\log y)}\Bigr)\biggr| \\
&&\qquad \ll \frac{\log^{2}y}{\log^{2}x} \sum_{0 \leq i < j \leq B} \sum
_{1 \leq k,l \leq M} \biggl|\E Z_{y}\biggl(2i+1+\frac{kE}{\log x}\biggr)
Z_{y}\biggl(2j+1+\frac{lE}{\log x}\biggr)\biggr|
\\
&&\qquad \ll \frac{\log^{2}y M^{2}}{\log^{2}x \log\log x} \sum_{0 \leq
i < j \leq B} \biggl(\frac{1}{|i-j|\log y} + (i+j)e^{-d\sqrt{\log y}} \biggr)
\\
&&\qquad \ll \frac{1}{(\log\log x)^{2}} \biggl(\frac{B\log B}{\log y} +
B^{3}e^{-d\sqrt{\log y}} \biggr).
\end{eqnarray*}

We noted above that, at the level of precision required in Section~\ref{sec6.2},
the correlation structure of $\{Z_{y}(t)\}_{t \in\mathcal{T}_{n}}$ is
the same for each $0 \leq n \leq B$. Thus our calculations concerning
$\sup_{t \in\mathcal{T}_{0}} Z_{y}(t)$ go through for $\sup_{t \in
\mathcal{T}_{n}} Z_{y}(t)$ as well, so that $\p(\sup_{t \in\mathcal
{T}_{n}} Z_{y}(t) \leq\sqrt{2(\log\log x - \log\log y)}) \leq
e^{-\Theta(\sqrt{\log\log\log x}/(\log\log x)^{2})}$ for each $0
\leq n \leq B$, and
\begin{eqnarray*}
&&\p\Bigl(\max_{0 \leq n \leq B} \sup_{t \in\mathcal
{T}_{n}} Z_{y}(t) \leq\sqrt{2(\log\log x - \log\log y)}\Bigr)\\
&&\qquad
\ll\frac{1}{(\log\log x)^{2}} \biggl(\frac{B\log B}{\log y} +
B^{3}e^{-d\sqrt{\log y}} \biggr)\\
&&\qquad\quad{} + e^{-\Theta((B+1)\sqrt{\log\log\log
x}/(\log\log x)^{2})}.
\end{eqnarray*}
The right-hand side is $O(e^{-\Theta(\sqrt{\log\log\log x})})$ if
we take $B=(\log\log x)^{2}$ and $y \geq\log x$.

For our application to $M(x)$, we need a version of the above
probability estimate in which $\max_{0 \leq n \leq B} \sup_{t \in
\mathcal{T}_{n}} Z_{y}(t)$ is replaced by
\[
\max_{0 \leq n \leq B} \sup_{t \in\mathcal{T}_{n}} \frac{\sum_{y
\leq p \leq x} f(p)\cos(t \log p)/p^{1/2+1/\log x}}{\sqrt{\sum_{y
\leq p \leq x} \cos^{2}(t \log p)/p^{1+2/\log x}}}
\]
with $f(p)$ independent Rademacher random variables. This can be
achieved using a multivariate central limit theorem, as explained in
Appendix~\ref{secB}, if we replace the upper bound $\sqrt{2(\log\log x - \log
\log y)}$ that we demand by $\sqrt{2(\log\log x - \log\log y) - 1}$.
In the application of the central limit theorem, we need $y$ to be at
least a certain power of $\log x$, say $y=\log^{8}x$. This choice is
also permissible for all of the preceding calculations.

Finally, note that for fixed $t \in\R$,
\begin{eqnarray*}
\E\biggl(\sum_{p < y} \frac{g_{p}\cos(t \log p)}{p^{1/2+1/\log x}}\biggr)^{2}
&=&
\E\biggl(\sum_{p < y} \frac{f(p)\cos(t \log p)}{p^{1/2+1/\log x}}\biggr)^{2} =
O(\log\log y)\\
&=& O(\log\log\log x)
\end{eqnarray*}
as $x \rightarrow\infty$, as in Section~\ref{sec6.1}. Applying Chebyshev's
inequality to this estimate,
\[
\p\biggl(\biggl|\sum_{p < y} \frac{g_{p}\cos(t \log p)}{p^{1/2+1/\log x}}\biggr| >
(\log\log\log x)^{3/4}\biggr) = O((\log\log\log x)^{-1/2}),
\]
also if the $g_{p}$ are replaced by Rademacher random variables
$f(p)$. These sums are independent of the sums over $y \leq p \leq x$,
so temporarily setting $d(x):=\inf_{1 \leq t \leq2(\log\log x)^{2}}
\sqrt{\E(\sum_{y \leq p \leq x} g_{p}\cos(t\log p)/p^{1/2+1/\log
x})^{2}}$ we find
\begin{eqnarray*}
&&
\p\biggl(\frac{1}{d(x)} \sup_{1 \leq t \leq2(\log\log x)^{2}} \sum_{p
\leq x} \frac{g_{p}\cos(t \log p)}{p^{1/2+1/\log x}} \\
&&\qquad\leq\sqrt
{2(\log\log x - \log\log y)} - \frac{(\log\log\log
x)^{3/4}}{d(x)}\biggr)
\end{eqnarray*}
is $O((\log\log\log x)^{-1/2})$. Corollary~\ref{co2} quickly follows since,
by the calculations in Section~\ref{sec6.1}, we have $d(x) = \sqrt{(\log
\log x
- \log\log y)/2 + O(1)}$.\qed

As noted in Appendix~\ref{secA}, the tail sum $\sum_{p > x} f(p)\cos(t \log
p)/p^{1/2+1/\log x}$ is almost surely convergent (and in fact it
converges in square mean) so that
\begin{eqnarray*}
\E\biggl(\sum_{p > x} \frac{f(p)\cos(t \log p)}{p^{1/2+1/\log x}}\biggr)^{2}
&\leq&\sum_{p > x} \frac{1}{p^{1+2/\log x}} \\
&=& O\biggl(\int_{x}^{\infty}
\frac{du}{u^{1+2/\log x} \log u}\biggr) \\
&=& O(1).
\end{eqnarray*}
Applying Chebyshev's inequality again, together with the Rademacher
version of our estimate for $Z_{y}(t)$, we have that
\begin{eqnarray*}
&&
\p\biggl(\sup_{1 \leq t \leq2(\log\log x)^{2}} \sum_{p} \frac{f(p)\cos
(t \log p)}{p^{1/2+1/\log x}}\\
&&\qquad \leq\log\log x - \log\log y - O(1) -
(\log\log\log x)^{3/4}\biggr)
\end{eqnarray*}
is $O((\log\log\log x)^{-1/2})$. Applying the first Borel--Cantelli
lemma at a lacunary set of points $x$, one quickly deduces that for any
fixed $A > 3$, there almost surely exists a sequence $(x_{k})$, tending
to infinity, with
\begin{eqnarray*}
&&
\sup_{1 \leq t \leq2(\log\log x_{k})^{2}} \sum_{p} \frac{f(p)\cos
(t \log p)}{p^{1/2+1/\log x_{k}}} - 2\log\log\log x_{k} \\
&&\qquad\geq\log
\log x_{k} - A\log\log\log x_{k}.
\end{eqnarray*}
By the argument in Appendix~\ref{secA} (and specifically by Supplementary
Lem\-ma~\ref{su1} from that appendix), this implies Corollary~\ref{co3} for $A
> 3$.

\section{\texorpdfstring{Refinement of Proposition \protect\ref{pr2} for the random multiplicative functions application}
{Refinement of Proposition 2 for the random multiplicative functions application}}\label{sec7}
As discussed at the end of Section~\ref{sec4}, Proposition~\ref{pr2}
may be refined in
that the product term can be replaced by any lower bound for
\[
\p\biggl(A_{j} \leq\frac{(1-\delta)(u-r_{j,m}(u+h))}{\sqrt{1 -
r_{j,m}^{2} - c_{j}d_{j}}}\ \forall1 \leq j \leq m-1 \biggr)
\]
for \textit{any} standard normal random variables $A_{j}$ satisfying
\[
\E A_{j}A_{k} \leq\frac{r_{j,k}-r_{j,m}r_{k,m}-c_{\min\{j,k\}
}d_{\max\{j,k\}}}{\sqrt
{(1-r_{j,m}^{2}-c_{j}d_{j})(1-r_{k,m}^{2}-c_{k}d_{k})}}.
\]
It will be convenient to write $U(j,k)$ for this upper bound on the
permissible correlations. By assumption about the numbers
$c_{j},d_{j}$, we always have \mbox{$U(j,k) \geq0$}.

In our application to random multiplicative functions, $U(j,k)$ is at least
\begin{eqnarray*}
&&
\frac{(r_{j,k}-1) + (1-r_{j,m})}{\sqrt
{(1-r_{j,m}^{2}-c_{j}d_{j})(1-r_{k,m}^{2}-c_{k}d_{k})}} \\ 
&&\qquad= \frac{-{\log}
|j-k|E + {\log}|j-m|E + O(1)}{\sqrt{({\log}|j-m|E+O(1))({\log}
|k-m|E+O(1))}}
\end{eqnarray*}
for $1 \leq j < k \leq m-1$. It seems sensible to consider intervals
$L^{i}/E < |m-j|$, $|m-k| \leq L^{i+1}/E$ (with $L \leq2$ a parameter to
be chosen) on which we see
\[
U(j,k) \geq1 - \frac{\log(|j-k|E)}{i\log L} + O\biggl(\frac{1}{i\log L}\biggr).
\]
In the random multiplicative functions example, on such an interval the
upper bound $(1-\delta)(u-r_{j,m}(u+h))/\sqrt{1 - r_{j,m}^{2} -
c_{j}d_{j}}$ that we demand for the $A_{j}$ is at least $(1+O(1/i\log
L))\sqrt{2i\log L}$. Thus, taking $A_{j}$ on distinct intervals to be
independent of one another (rather than \textit{all} $A_{j}$ necessarily
being independent), we can replace the product in Proposition~\ref{pr2} by
\[
\mathop{\prod_{i=0,}}_{L^{i} \geq E/2}^{[\log(Em)/\log L]} \p\biggl(A_{j} \leq
\biggl(1 - \frac{c}{i\log L}\biggr)\sqrt{2i\log L}\ \forall L^{i}/E < |m-j| \leq
L^{i+1}/E\biggr),
\]
where $c$ is an absolute constant and $A_{j}$ are any standard normal
random variables whose correlations are bounded as described.

The crucial point is that on each interval, and up to the ``big Oh''
term, the bound on $U(j,k)$ corresponds to a stationary correlation
structure that we can hope to understand. Indeed, it is essentially a
re-scaled version of the original correlation structure of our random
multiplicative functions process.

Using these ideas, we shall establish the following result. In its
statement we include a superscript $x$ to explicitly record that
$Z_{y}(t)=Z_{y}^{x}(t)$ depends on $x$, and we remind the reader that
we had $y=\log^{8}x$.
%
%
\begin{proposition}\label{pr3}
If $E$ is a sufficiently large constant, then the following is true.
Let $\{Z_{y}(t)\}_{t \in\mathcal{T}_{0}} = \{Z_{y}^{x}(t)\}_{t \in
\mathcal{T}_{0}}$ be the Gaussian process described in Section~\ref{sec6.1},
for such a choice of $E$. Let $\varepsilon(x)$ be any function tending to
zero as $x \rightarrow\infty$. Then for some sequence of $x$, tending
to infinity, we have
\[
\p\biggl(\sup_{t \in\mathcal{T}_{0}} Z_{y}(t) > \sqrt{2(\log\log x -
\log\log y)}\biggr) \geq\frac{\varepsilon(x) \sqrt{\log E}}{E(\log\log
x)^{3/2}}.
\]
\end{proposition}

Recall from Section~\ref{sec6.1} that
\begin{eqnarray}
\E Z_{y}^{x}\biggl(1+\frac{jE}{\log x}\biggr)Z_{y}^{x}\biggl(1+\frac{kE}{\log x}\biggr) = 1 -
\frac{\log(|j-k|E) + O(1)}{\log\log x - \log\log y},
\nonumber\\
&&\eqntext{\displaystyle 1 \leq j,k
\leq\frac{\log x}{KE\log y}, j \neq k,}
\end{eqnarray}
where $K$ is an absolute constant in the definition of $\mathcal
{T}_{0}$. Let us fix a large absolute constant $C \in\N$, and set $L
= 1+1/KC^{3}$. When $L^{i}$ is large enough, we can choose $x(i) \in\R
$ such that
\[
\sqrt{2\bigl(\log\log x(i) - \log\log y(i)\bigr)} = \biggl(1 -
\frac{c}{i \log L}\biggr) \sqrt{2i\log L}.
\]
Here we wrote $y(i)=y(x(i))=\log^{8}x(i)$. Then we will have
\begin{eqnarray*}
\E Z_{y(i)}^{x(i)}\biggl(1+\frac{jCE}{\log
x(i)}\biggr)Z_{y(i)}^{x(i)}\biggl(1+\frac {kCE}{\log x(i)}\biggr) &=& 1 -
\frac{\log(|j-k|E) + \log C + O(1)}{(1-c/(i\log L))^{2} i \log L} \\
&\leq&
U(j,k),
\end{eqnarray*}
where $U(j,k)$ denotes the bound for interval $i$. This only makes
sense if $jC, kC \leq\log x(i) / KE\log y(i)$, but that will hold, for
example,\vadjust{\goodbreak} if $j, k \leq L^{i}/\break KEC^{2}$. Thus if $i$ is sufficiently
large that
\[
\biggl[\frac{L^{i+1}}{E}\biggr] - \biggl[\frac{L^{i}}{E}\biggr] \leq\biggl[\frac
{L^{i}}{KEC^{2}}\biggr],
\]
we can say that $\p(A_{j} \leq(1 - \frac{c}{i\log L})\sqrt{2i\log
L}\ \forall L^{i}/E < |m-j| \leq L^{i+1}/E)$ is at least
\[
\p\Bigl(\sup_{t \in\mathcal{T}_{0}} Z_{y(i)}^{x(i)}(t) \leq\sqrt
{2\bigl(\log\log x(i) - \log\log y(i)\bigr)}\Bigr).
\]

Notice that, for our fixed choice of $L$, the various requirements for
$i$ to be ``sufficiently large'' will all be satisfied if $i \geq
i_{E}+D$, where $i_{E}$ is least for which $L^{i} \geq E/2$ and
$D=D(L)$ is a constant. Thus the product term in Proposition~\ref{pr2}
may be
replaced by
\begin{eqnarray*}
&& \prod_{j=1}^{[L^{D}]} \Phi\biggl( \biggl(1+O\biggl(\frac{1}{\log jE}\biggr)\biggr) \sqrt{2\log
jE} \biggr) \\
&&\qquad\hspace*{0pt}{}
\times\prod_{i=i_{E}+D}^{[{\log(Em)}/{\log L}]} \p\Bigl(\sup_{t
\in\mathcal{T}_{0}} Z_{y(i)}^{x(i)}(t) \leq\sqrt{2\bigl(\log\log x(i) -
\log\log y(i)\bigr)}\Bigr).
\end{eqnarray*}
We also note that, obviously, $x(i)$ tends to infinity with $i$.

Now suppose that the proposition failed, so for all sufficiently large
$x$ the tail probability was smaller than required. Then for all $i$
from some point onward we would have
\begin{eqnarray*}
\p\Bigl(\sup_{t \in\mathcal{T}_{0}} Z_{y(i)}^{x(i)}(t) \leq\sqrt
{2\bigl(\log\log x(i) - \log\log y(i)\bigr)}\Bigr)
&\geq& 1 - \frac{1}{(\log\log
x(i))^{3/2}} \\
& \geq& 1 - O\biggl(\frac{1}{(i\log L)^{3/2}}\biggr),
\end{eqnarray*}
so [since $\prod_{i=2}^{\infty}(1-1/i^{3/2})$ is convergent] the
product term in Proposition~\ref{pr2} could be replaced by a positive constant.
But then the argument of Section~\ref{sec6.2} would supply that
\[
\p\Bigl(\sup_{t \in\mathcal{T}_{0}} Z_{y}(t) > \sqrt{2(\log\log x -
\log\log y)}\Bigr) \gg\frac{\sqrt{\log E}}{E(\log\log x)^{3/2}},
\]
which is a contradiction for $x$ sufficiently large.\qed\vspace*{8pt}
%

Armed with Proposition~\ref{pr3}, we can repeat the argument of Section
\ref{sec6.3}
with $E$ chosen to be a large constant (rather than $\sqrt{\log\log
x}$), and $B$ then chosen as $(\log\log x)^{3/2} \log\log\log x$,
say [rather than $(\log\log x)^{2}$]. The reader should note that
there is a subtlety involved, as this requires lower bounds for
\[
\p\Bigl(\sup_{t \in\mathcal{T}_{n}} Z_{y}(t) > \sqrt{2(\log\log x -
\log\log y)}\Bigr),\qquad 0 \leq n \leq B,\vadjust{\goodbreak}
\]
while Proposition~\ref{pr3} concerns $\sup_{t \in\mathcal{T}_{0}} Z_{y}(t)$
only. However, modifying the choice of $E$ and $K$ by some
multiplicative constants in the definition of $\mathcal{T}_{n}$, $n
\neq0$, so that $E$ is larger but $EK$ remains the same, we can
arrange using Comparison Inequality~\ref{comp2} that
\begin{eqnarray*}
&&\p\Bigl(\sup_{t \in\mathcal{T}_{n}} Z_{y}(t) > \sqrt{2(\log\log x -
\log\log y)}\Bigr)\\
&&\qquad \geq\p\Bigl(\sup_{t \in\mathcal{T}_{0}} Z_{y}(t) > \sqrt
{2(\log\log x - \log\log y)}\Bigr).
\end{eqnarray*}
We also only have probability bounds for a sequence of $x$ tending to
infinity, rather than all $x$, but we do not require that in Section
\ref{sec6.3}. Corollary~\ref{co3} follows from these considerations.

%
\begin{appendix}\label{app}
\section{Random multiplicative functions and Rademacher processes}\label{secA}
In this Appendix we sketch the connection between the sum $M(x) = \sum
_{n \leq x} f(n)$ of a random multiplicative function (as defined in
the \hyperref[sec1]{Introduction}) and a certain Rademacher random process. The argument
we give is essentially that of Hal\'{a}sz~\cite{hal}.

In view of Wintner's~\cite{wint} result that for each $\varepsilon> 0$,
$M(x)=O(x^{1/2+\varepsilon})$ almost surely, we know that the Dirichlet series
\[
F(s):= \sum_{n=1}^{\infty} \frac{f(n)}{n^{s}}
\]
is almost surely convergent in the half plane $\Re(s) > 1/2$, and then
satisfies
\[
F(s) = s \int_{1}^{\infty} \frac{M(z)}{z^{s+1}} \,dz.
\]
On the other hand, writing $\zeta(s):=\sum_{n}1/n^{s}, \Re(s) > 1$
for the Riemann zeta function, we have the Euler product identity
\begin{eqnarray*}
F(s) &=& \prod_{p} \biggl(1+\frac{f(p)}{p^{s}}\biggr) \\
& = & e^{\sum_{p} f(p)/p^{s}
- \sum_{p} 1/2p^{2s} + \sum_{k \geq3} \sum_{p} (-1)^{k+1}
f(p)^{k}/k p^{ks}} \\
& = & e^{\sum_{p} f(p)/p^{s} - \log\zeta(2s)/2 + \sum_{k \geq2}
\sum_{p} 1/2k p^{2ks} + \sum_{k \geq3} \sum_{p} (-1)^{k+1}
f(p)^{k}/k p^{ks}}.
\end{eqnarray*}
This is certainly valid when $\Re(s) > 1$, and almost surely extends
to $\Re(s) > 1/2$ in view of Kolmogorov's three series theorem and the
identity theorem of complex analysis. [The three series theorem implies
that $\sum_{p} f(p)/p^{s}$ converges almost surely when $\Re(s) >
1/2$. We then use the standard fact, proved using partial summation,
that such a Dirichlet series is a holomorphic function strictly to the
right of its abscissa of convergence.]\vadjust{\goodbreak}

Thus in the domain $1/2 < \sigma< 1, 1 \leq t \leq2$, say, we almost
surely have
\begin{eqnarray*}
\frac{e^{\sum_{p} f(p) \cos(t \log p)/p^{\sigma}}}{t} &\ll&\int
_{1}^{\infty} \frac{|M(z)|}{z^{\sigma+1}} \,dz \\
&\leq&\sup_{z \geq1}
\frac{|M(z)|}{\sqrt{z(\sigma-1/2)}} + \sup_{z \geq z_{0}} \frac
{|M(z)|}{\sqrt{z} (\sigma-1/2)},
\end{eqnarray*}
where the second inequality follows by splitting the integral at
$z_{0}:=\break e^{1/\sqrt{\sigma-1/2}}$. Taking $\sigma=1/2+1/\log x$,
where $x \geq2$ is a parameter, we find that
\begin{eqnarray*}
&&
e^{\sum_{p} f(p) \cos(t \log p)/p^{1/2+1/\log x}} \\
&&\qquad\ll\sqrt{\log x}
\sup_{z \geq1} \frac{|M(z)|}{\sqrt{z}} + \log x \sup_{z \geq
e^{\sqrt{\log x}}} \frac{|M(z)|}{\sqrt{z}},\qquad 1 \leq t \leq2.
\end{eqnarray*}

For the proof of Corollary~\ref{co3}, we need a version of the preceding
inequality that is valid for a larger range of $t$. Using the estimate
$|{\log\zeta}(\sigma+it)| \leq{\log\log}|t| + O(1), \sigma\geq1,
|t| \geq2$, which is contained in, for example, Montgomery and Vaughan
\cite{mv}, Theorem 6.7, we can say that for $t \geq1$,
\begin{eqnarray*}
&&
e^{\sum_{p} f(p)\cos(t \log p)/p^{1/2+1/\log x} - \log t - \log\log
(t+2)/2} \\
&&\qquad\ll\sqrt{\log x} \sup_{z \geq1} \frac{|M(z)|}{\sqrt{z}}
+ \log x \sup_{z \geq e^{\sqrt{\log x}}} \frac{|M(z)|}{\sqrt{z}}.
\end{eqnarray*}
This immediately implies the following result.
%
%
\begin{supplem}\label{su1}
Let $g(z)$ be a decreasing function. If, with positive probability, we
have $M(z) = O(\sqrt{z}g(z))$ as $z \rightarrow\infty$, then with
positive probability we have
\begin{eqnarray*}
&&
\sup_{t \geq1} e^{\sum_{p} f(p)\cos(t \log p)/p^{1/2+1/\log x} -
\log t - \log\log(t+2)/2} \\
&&\qquad= O\bigl(g(1)\sqrt{\log x} + g\bigl(e^{\sqrt{\log
x}}\bigr)\log x\bigr)
\end{eqnarray*}
for all $x \geq2$.
\end{supplem}

Since Hal\'{a}sz's paper~\cite{hal} seems to be difficult to get hold
of, it is perhaps worthwhile to briefly discuss Hal\'{a}sz's own use of
the foregoing argument. He shows that there almost surely exist
sequences of real numbers $x_{k}$, tending to infinity, and of sets
$S_{k} \subseteq[1,2]$, of measure $> 1/\log x_{k}$ and of sets $B_{k}
\subseteq[1,2]$, of measure $\leq1/\log x_{k}$, such that
\begin{eqnarray*}
&&\sum_{p \leq x_{k}} f(p) \frac{\cos(t \log p)}{\sqrt{p}} \\
&&\qquad\geq\log
\log x_{k} - \sqrt{29\log\log x_{k} \log\log\log x_{k}}\qquad \forall t
\in S_{k},
\end{eqnarray*}
and
\begin{eqnarray*}
&&\sum_{p \leq x_{k}} f(p) \frac{\cos(t \log p)}{\sqrt{p}} - \sum
_{p} f(p) \frac{\cos(t \log p)}{p^{1/2+1/\log x_{k}}} \\
&&\qquad= O\bigl(\sqrt{\log
\log x_{k}}\bigr)\qquad \forall t \in[1,2] \setminus B_{k}.
\end{eqnarray*}
In particular, there almost surely exists a sequence $x_{k}$ such that
\begin{eqnarray*}
&&\sup_{t \in[1,2]} \sum_{p} f(p) \frac{\cos(t \log
p)}{p^{1/2+1/\log x_{k}}} \\
&&\qquad\geq\log\log x_{k} - \sqrt{29\log\log
x_{k} \log\log\log x_{k}} - O\bigl(\sqrt{\log\log x_{k}}\bigr),
\end{eqnarray*}
which is enough to imply the omega result for $M(x)$ attributed to
Hal\'{a}sz in the \hyperref[sec1]{Introduction}.

Very roughly, Hal\'{a}sz~\cite{hal} investigates the process $\sum_{p
\leq x} f(p) \frac{\cos(t \log p)}{\sqrt{p}}$, $t \in[1,2]$, by
estimating moments of the counting function
\[
\int_{1}^{2} \mathbf{1}_{\sum_{p \leq x} f(p) \cos(t \log p)/\sqrt
{p} \geq M} \,dt,
\]
where $M$ is a parameter. However, the details are rather complicated,
as it is actually necessary to split the sum over $p$ into several
ranges, and then reduce the range of integration to progressively
smaller random subsets of $[1,2]$. This splitting is, in a sense, quite
natural, as the parts of the sum taken over large primes are less
correlated at nearby values of $t$ (see Section~\ref{sec6.1}). On the other
hand, the splitting causes an accumulation of error terms in the
analysis, one from each range of summation. The iterative approach is
also highly reliant on being presented with the process as a random sum
over~$p$, whereas [at least if the $f(p)$ were independent Gaussians]
one might just as well be given a description of the process only in
terms of its covariance structure.\looseness=-1

\section{A multivariate central limit theorem}\label{secB}
In this Appendix we discuss a multivariate central limit theorem of
Reinert and R\"{o}llin~\cite{rr}. We view this as a ``universality
result,'' which sometimes lets us transfer conclusions about suprema of
Gaussian processes to conclusions about the suprema of corresponding
Rademacher processes. Reinert and R\"{o}llin's~\cite{rr} approach is
based on Stein's method of exchangeable pairs.

Suppose that $\mathcal{T}$ is a finite set, and that $\alpha_{i}(t)
\in\R$ for $1 \leq i \leq n$ and $t \in\mathcal{T}$. Suppose also
that $(\varepsilon_{i})_{i=1}^{n}$ is a sequence of independent
Rademacher random variables and that $(g_{i})_{i=1}^{n}$ is a sequence
of independent standard normal random variables. We wish to approximate
the (joint) distribution of $\{X_{t}\}_{t \in\mathcal{T}}$ by that of
$\{Y_{t}\}_{t \in\mathcal{T}}$, where
\[
X_{t}:= \sum_{i=1}^{n} \alpha_{i}(t) \varepsilon_{i},\qquad Y_{t}:= \sum
_{i=1}^{n} \alpha_{i}(t) g_{i}.
\]

In the usual way, we construct random variables $X_{t}'$ so
$((X_{t})_{t \in\mathcal{T}},(X_{t}')_{t \in\mathcal{T}})$ is an
exchangeable pair of vectors [i.e., so that the law of this tuple
is the same as the law of $((X_{t}')_{t \in\mathcal{T}},(X_{t})_{t
\in\mathcal{T}})$]. Let $I$ be a random variable having the discrete
uniform distribution on $\{1,2,\ldots,n\}$, independently of everything
else and let $(\varepsilon_{i}')_{i=1}^{n}$ be an independent copy of
$(\varepsilon_{i})_{i=1}^{n}$. We define $X_{t}'$ as follows: conditional
on the event $\{I=i\}$, set
\[
X_{t}'=X_{t} - \alpha_{i}(t)\varepsilon_{i} + \alpha_{i}(t)\varepsilon
_{i}',\qquad t \in\mathcal{T}.
\]
The reader may check that the exchangeability property does then hold,
together with the following regression property:
\[
\E\bigl(X_{t}'-X_{t}|(X_{s})_{s \in\mathcal{T}}\bigr) = -\frac{1}{n}X_{t}.
\]

With a view to applying Theorem 2.1 of Reinert and R\"{o}llin \cite
{rr}, we calculate two further quantities:
\begin{eqnarray*}
\E\bigl((X_{t}'-X_{t})(X_{s}'-X_{s}) | (X_{u})_{u \in\mathcal{T}}\bigr) & = &
\frac{1}{n} \sum_{i=1}^{n} \alpha_{i}(t)\alpha_{i}(s) \E\bigl((\varepsilon
_{i}'-\varepsilon_{i})^{2}|(X_{u})_{u \in\mathcal{T}}\bigr) \\
& = & \frac{2}{n} \sum_{i=1}^{n} \alpha_{i}(t)\alpha_{i}(s);
\\
\E|(X_{t}'-X_{t})(X_{s}'-X_{s})(X_{u}'-X_{u})| & = & \frac{1}{n} \sum
_{i=1}^{n} |\alpha_{i}(t)\alpha_{i}(s)\alpha_{i}(u)| \E|\varepsilon
_{i}'-\varepsilon_{i}|^{3} \\
& = & \frac{4}{n} \sum_{i=1}^{n} |\alpha_{i}(t)\alpha_{i}(s)\alpha
_{i}(u)|.
\end{eqnarray*}
The reader should notice that, while we did not use the fact that the
$\varepsilon_{i}$ are Rademacher random variables up until this point, in
the first calculation it allows us to conclude that the left-hand side
is deterministic. This means that one of the error terms in Reinert and
R\"{o}llin's~\cite{rr} theorem is identically zero; indeed, if $h\dvtx\R
^{\#\mathcal{T}} \rightarrow\R$ is a three times differentiable
function, and if the covariance matrix of $(X_{t})_{t \in\mathcal
{T}}$ is nonsingular, their theorem implies that
\begin{eqnarray*}
&&
|\E h((X_{t})_{t \in\mathcal{T}}) - \E h((Y_{t})_{t \in\mathcal
{T}})| \\
&&\qquad\leq\frac{1}{3} \sup_{s,t,u \in\mathcal{T}, \tilde{x} \in
\R^{\#\mathcal{T}}} \biggl|\frac{\partial^{3} h(\tilde{x})}{\partial
x_{s} \,\partial x_{t} \,\partial x_{u}} \biggr| \sum_{s,t,u \in\mathcal{T}}
\sum_{i=1}^{n} |\alpha_{i}(s) \alpha_{i}(t) \alpha_{i}(u)|.
\end{eqnarray*}
The condition that the covariance matrix should be nonsingular is
evidently unnecessary here (at least if $h$ is bounded, say), since we
can ensure this by introducing $\#\mathcal{T}$ dummy random variables
whose coefficients $\alpha_{i}(t)$ have absolute value at most~$\delta
$, and then let $\delta\rightarrow0$.

Specializing to our random multiplicative functions application, we
would like to choose $h$ to be the indicator function of a box in $\R
^{\#\mathcal{T}}$, but this would not satisfy the three times
differentiability condition. Reinert and R\"{o}llin devote a section of
their paper~\cite{rr} to this ``unsmoothing'' problem, but the results
they obtain are rather involved, and in this case we can easily
overcome the difficulty directly. Let $s\dvtx\R\rightarrow[0,1]$ be a
three times differentiable function satisfying
\[
s(z) = \cases{
1, &\quad if $z \leq\sqrt{2(\log\log x - \log\log y) -1}$, \vspace*{2pt}\cr
0, &\quad if $z \geq\sqrt{2(\log\log x - \log\log y)}$.}
\]
The interval\vspace*{1pt} on which $s(z)$ must drop from $1$ to $0$ has length
$\Theta(1/\sqrt{\log\log x})$, so we
can find such $s$ with
derivatives satisfying $|s^{(r)}(z)| = O((\log\log x)^{r/2})$, $0 \leq
r \leq3$, $z \in\R$. Setting $h((x_{t})_{t \in\mathcal{T}}) =
\prod_{t \in\mathcal{T}} s(x_{t})$, we conclude that
\begin{eqnarray*}
&&
\p\Bigl(\max_{t \in\mathcal{T}} X_{t} \leq\sqrt{2(\log\log x - \log
\log y) -1}\Bigr)
\\
&&\qquad\leq\p\Bigl(\max_{t \in\mathcal{T}} Y_{t} \leq\sqrt{2(\log\log x -
\log\log y)}\Bigr)\\
&&\qquad\quad{} + O\Biggl((\log\log x)^{3/2} (\#\mathcal{T})^{3} \sum
_{i=1}^{n} {\max_{t \in\mathcal{T}} }|\alpha_{i}(t)|^{3}\Biggr).
\end{eqnarray*}
The reader may check that in the random multiplicative functions case,
the error term on the right-hand side has order at most
\[
(\#\mathcal{T})^{3} \sum_{y \leq p \leq x} \frac{1}{p^{3/2}} \ll
\frac{(\#\mathcal{T})^{3}}{\sqrt{y} \log y}.
\]
We have $\#\mathcal{T} = (B+1)M \ll(\log\log x)^{2} \log x$, so this
is $o(1)$ as $x \rightarrow\infty$ provided that $y$ is at least
$\log^{8} x$, say. The multivariate central limit theorem has supplied
an extremely good bound, presumably because any individual $\varepsilon
_{p}$ (or $g_{p}$) has a very tiny impact on the random multiplicative
function processes.
\end{appendix}

\section*{Acknowledgments}

The author would like to thank his Ph.D. supervisor, Ben Green, for
introducing him to Hal\'{a}sz's work on random multiplicative
functions, and for reading a draft of this paper. He would also like to
thank Nathana\"{e}l Berestycki and Richard Nickl for discussions on
these topics and the anonymous referee for his or her comments.



\printaddresses


\begin{thebibliography}{22}

\bibitem{burnmich}
\begin{barticle}[mr]
\bauthor{\bsnm{Burnecki},~\bfnm{K.}\binits{K.}} \AND
  \bauthor{\bsnm{Michna},~\bfnm{Z.}\binits{Z.}}
(\byear{2002}).
\btitle{Simulation of {P}ickands constants}.
\bjournal{Probab. Math. Statist.}
\bvolume{22}
\bpages{193--199}.
\bid{issn={0208-4147}, mr={1944151}}
\bptok{imsref}%
\end{barticle}
\endbibitem

\bibitem{chatsound}
\begin{barticle}[auto:STB|2012/05/30|10:51:56]
\bauthor{\bsnm{Chatterjee},~\bfnm{S.}\binits{S.}} \AND
  \bauthor{\bsnm{Soundararajan},~\bfnm{K.}\binits{K.}}
(\byear{2012}).
\btitle{Random multiplicative functions in short intervals}.
\bjournal{Int. Math. Res. Not.}
\bvolume{3}
\bpages{479--492}.
\bptok{imsref}%
\end{barticle}\vadjust{\goodbreak}
\endbibitem

\bibitem{dk}
\begin{barticle}[mr]
\bauthor{\bsnm{D{\c{e}}bicki},~\bfnm{Krzysztof}\binits{K.}} \AND
  \bauthor{\bsnm{Kisowski},~\bfnm{Pawe{\l}}\binits{P.}}
(\byear{2008}).
\btitle{A note on upper estimates for {P}ickands constants}.
\bjournal{Statist. Probab. Lett.}
\bvolume{78}
\bpages{2046--2051}.
\bid{doi={10.1016/j.spl.2008.01.071}, issn={0167-7152}, mr={2458013}}
\bptok{imsref}%
\end{barticle}
\endbibitem

\bibitem{dmr}
\begin{barticle}[mr]
\bauthor{\bsnm{D\c{e}bicki},~\bfnm{Krzysztof}\binits{K.}},
  \bauthor{\bsnm{Michna},~\bfnm{Zbigniew}\binits{Z.}} \AND
  \bauthor{\bsnm{Rolski},~\bfnm{Tomasz}\binits{T.}}
(\byear{2003}).
\btitle{Simulation of the asymptotic constant in some fluid models}.
\bjournal{Stoch. Models}
\bvolume{19}
\bpages{407--423}.
\bid{doi={10.1081/STM-120023567}, issn={1532-6349}, mr={1993949}}
\bptok{imsref}%
\end{barticle}
\endbibitem

\bibitem{gs}
\begin{bbook}[mr]
\bauthor{\bsnm{Grimmett},~\bfnm{Geoffrey~R.}\binits{G.~R.}} \AND
  \bauthor{\bsnm{Stirzaker},~\bfnm{David~R.}\binits{D.~R.}}
(\byear{2001}).
\btitle{Probability and Random Processes}, \bedition{3rd} ed.
\bpublisher{Oxford Univ. Press}, \baddress{New York}.
\bid{mr={2059709}}
\bptok{imsref}%
\end{bbook}
\endbibitem

\bibitem{hal}
\begin{bincollection}[mr]
\bauthor{\bsnm{Hal{\'a}sz},~\bfnm{G.}\binits{G.}}
(\byear{1983}).
\btitle{On random multiplicative functions}.
In \bbooktitle{Hubert {D}elange Colloquium ({O}rsay, 1982)}.
\bseries{Publications Math\'{e}matiques d'Orsay}
\bvolume{83}
\bpages{74--96}.
\bpublisher{Univ. Paris XI}, \baddress{Orsay}.
\bid{mr={0728404}}
\bptok{imsref}%
\end{bincollection}
\endbibitem

\bibitem{harper}
\begin{bmisc}[auto:STB|2012/05/30|10:51:56]
\bauthor{\bsnm{Harper},~\bfnm{A.~J.}\binits{A.~J.}}
(\byear{2012}).
\bhowpublished{On the limit distributions of some sums of a random
  multiplicative function. \textit{J. Reine Angew. Math.} To appear.}
\bptok{imsref}%
\end{bmisc}
\endbibitem

\bibitem{ho2}
\begin{barticle}[mr]
\bauthor{\bsnm{Hough},~\bfnm{Bob}\binits{B.}}
(\byear{2011}).
\btitle{Summation of a random multiplicative function on numbers having few
  prime factors}.
\bjournal{Math. Proc. Cambridge Philos. Soc.}
\bvolume{150}
\bpages{193--214}.
\bid{doi={10.1017/S0305004110000514}, issn={0305-0041}, mr={2770059}}
\bptok{imsref}%
\end{barticle}
\endbibitem

\bibitem{ltw}
\begin{barticle}[auto:STB|2012/05/30|10:51:56]
\bauthor{\bsnm{Lau},~\bfnm{Y.~K.}\binits{Y.~K.}},
  \bauthor{\bsnm{Tenenbaum},~\bfnm{G.}\binits{G.}} \AND
  \bauthor{\bsnm{Wu},~\bfnm{J.}\binits{J.}}
(\byear{2013}).
\btitle{On mean values of random multiplicative functions}.
\bjournal{Proc. Amer. Math. Soc.}
\bvolume{141}
\bpages{409--420}.
\bptok{imsref}%
\end{barticle}
\endbibitem

\bibitem{llr}
\begin{bbook}[mr]
\bauthor{\bsnm{Leadbetter},~\bfnm{M.~R.}\binits{M.~R.}},
  \bauthor{\bsnm{Lindgren},~\bfnm{Georg}\binits{G.}} \AND
  \bauthor{\bsnm{Rootz{\'e}n},~\bfnm{Holger}\binits{H.}}
(\byear{1983}).
\btitle{Extremes and Related Properties of Random Sequences and Processes}.
\bpublisher{Springer}, \baddress{New York}.
\bid{mr={0691492}}
\bptok{imsref}%
\end{bbook}
\endbibitem

\bibitem{lishao}
\begin{barticle}[mr]
\bauthor{\bsnm{Li},~\bfnm{Wenbo~V.}\binits{W.~V.}} \AND
  \bauthor{\bsnm{Shao},~\bfnm{Qi-Man}\binits{Q.-M.}}
(\byear{2002}).
\btitle{A normal comparison inequality and its applications}.
\bjournal{Probab. Theory Related Fields}
\bvolume{122}
\bpages{494--508}.
\bid{doi={10.1007/s004400100176}, issn={0178-8051}, mr={1902188}}
\bptok{imsref}%
\end{barticle}
\endbibitem

\bibitem{lif}
\begin{bbook}[mr]
\bauthor{\bsnm{Lifshits},~\bfnm{M.~A.}\binits{M.~A.}}
(\byear{1995}).
\btitle{Gaussian Random Functions}.
\bseries{Mathematics and Its Applications}
\bvolume{322}.
\bpublisher{Kluwer Academic}, \baddress{Dordrecht}.
\bid{mr={1472736}}
\bptok{imsref}%
\end{bbook}
\endbibitem

\bibitem{michna}
\begin{bmisc}[auto:STB|2012/05/30|10:51:56]
\bauthor{\bsnm{Michna},~\bfnm{Z.}\binits{Z.}}
(\byear{2009}).
\bhowpublished{Remarks on Pickands theorem. Preprint.}
\bptok{imsref}%
\end{bmisc}
\endbibitem

\bibitem{mv}
\begin{bbook}[mr]
\bauthor{\bsnm{Montgomery},~\bfnm{Hugh~L.}\binits{H.~L.}} \AND
  \bauthor{\bsnm{Vaughan},~\bfnm{Robert~C.}\binits{R.~C.}}
(\byear{2007}).
\btitle{Multiplicative Number Theory. {I}. {C}lassical Theory}.
\bseries{Cambridge Studies in Advanced Mathematics}
\bvolume{97}.
\bpublisher{Cambridge Univ. Press}, \baddress{Cambridge}.
\bid{mr={2378655}}
\bptok{imsref}%
\end{bbook}
\endbibitem

\bibitem{pickands1}
\begin{barticle}[mr]
\bauthor{\bsnm{Pickands},~\bfnm{James}\binits{J.} \bsuffix{III}}
(\byear{1969}).
\btitle{Upcrossing probabilities for stationary {G}aussian processes}.
\bjournal{Trans. Amer. Math. Soc.}
\bvolume{145}
\bpages{51--73}.
\bid{issn={0002-9947}, mr={0250367}}
\bptok{imsref}%
\end{barticle}
\endbibitem

\bibitem{pickands2}
\begin{barticle}[mr]
\bauthor{\bsnm{Pickands},~\bfnm{James}\binits{J.} \bsuffix{III}}
(\byear{1969}).
\btitle{Asymptotic properties of the maximum in a stationary {G}aussian
  process.}
\bjournal{Trans. Amer. Math. Soc.}
\bvolume{145}
\bpages{75--86}.
\bid{issn={0002-9947}, mr={0250368}}
\bptok{imsref}%
\end{barticle}
\endbibitem

\bibitem{pit}
\begin{bbook}[mr]
\bauthor{\bsnm{Piterbarg},~\bfnm{Vladimir~I.}\binits{V.~I.}}
(\byear{1996}).
\btitle{Asymptotic Methods in the Theory of {G}aussian Processes and Fields}.
\bseries{Translations of Mathematical Monographs}
\bvolume{148}.
\bpublisher{Amer. Math. Soc.}, \baddress{Providence, RI}.
\bid{mr={1361884}}
\bptok{imsref}%
\end{bbook}
\endbibitem

\bibitem{plac}
\begin{barticle}[mr]
\bauthor{\bsnm{Plackett},~\bfnm{R.~L.}\binits{R.~L.}}
(\byear{1954}).
\btitle{A reduction formula for normal multivariate integrals}.
\bjournal{Biometrika}
\bvolume{41}
\bpages{351--360}.
\bid{issn={0006-3444}, mr={0065047}}
\bptok{imsref}%
\end{barticle}
\endbibitem

\bibitem{rr}
\begin{barticle}[mr]
\bauthor{\bsnm{Reinert},~\bfnm{Gesine}\binits{G.}} \AND
  \bauthor{\bsnm{R{\"o}llin},~\bfnm{Adrian}\binits{A.}}
(\byear{2009}).
\btitle{Multivariate normal approximation with {S}tein's method of exchangeable
  pairs under a general linearity condition}.
\bjournal{Ann. Probab.}
\bvolume{37}
\bpages{2150--2173}.
\bid{doi={10.1214/09-AOP467}, issn={0091-1798}, mr={2573554}}
\bptok{imsref}%
\end{barticle}
\endbibitem

\bibitem{shao}
\begin{barticle}[mr]
\bauthor{\bsnm{Shao},~\bfnm{Qi-Man}\binits{Q.-M.}}
(\byear{1996}).
\btitle{Bounds and estimators of a basic constant in extreme value theory of
  {G}aussian processes}.
\bjournal{Statist. Sinica}
\bvolume{6}
\bpages{245--257}.
\bid{issn={1017-0405}, mr={1379060}}
\bptok{imsref}%
\end{barticle}
\endbibitem

\bibitem{wint}
\begin{barticle}[mr]
\bauthor{\bsnm{Wintner},~\bfnm{Aurel}\binits{A.}}
(\byear{1944}).
\btitle{Random factorizations and {R}iemann's hypothesis}.
\bjournal{Duke Math. J.}
\bvolume{11}
\bpages{267--275}.
\bid{issn={0012-7094}, mr={0010160}}
\bptok{imsref}%
\end{barticle}
\endbibitem

\end{thebibliography}
\end{document}